\documentclass[11pt,a4paper,reqno]{amsart}

\usepackage[T1]{fontenc}
\usepackage{lmodern}
\usepackage{amsmath,amsthm,amssymb}
\usepackage{mathtools}
\usepackage{mathrsfs}
\IfFormatAtLeastTF{2026-06-01}{}{\usepackage{aliascnt}}
\usepackage{enumitem}
\usepackage{microtype}
\AtBeginDocument{\begingroup\footnotesize\fontencoding{OT1}\selectfont\endgroup}
\usepackage[a4paper,textwidth=145mm,textheight=632pt,centering]{geometry}
\usepackage{xcolor}

\definecolor{mathCarmine}{HTML}{7A1F3D} % Deep carmine / berry for internal math refs
\definecolor{mathForest}{HTML}{1B5E3C}  % Deep forest green for literature citations
\definecolor{mathSapphire}{HTML}{0B4F8A}% Royal sapphire for web links and repositories

\usepackage{hyperref}
\usepackage{amsrefs}
\usepackage[nameinlink,capitalise]{cleveref}
\hypersetup{
    colorlinks = true,
    linkcolor  = mathCarmine,
    citecolor  = mathForest,
    urlcolor   = mathSapphire,
    bookmarksnumbered=true,
    pdftitle={C(m,alpha)-rectifiability in Carnot groups},
    pdfauthor={Kennedy Obinna Idu, Francesco Paolo Maiale},
    pdfsubject={Higher-order rectifiability in Carnot groups},
    pdfkeywords={higher-order rectifiability, Carnot groups, Whitney extension}
}

\numberwithin{equation}{section}

\makeatletter\def\bibtex@style{amsrndoi}\makeatother
\renewcommand{\PrintDOI}[1]{{\urlstyle{same}\url{https://doi.org/#1}}}

\newcommand{\eps}{\varepsilon}
\newcommand{\N}{\mathbb{N}}
\newcommand{\R}{\mathbb{R}}

\newcommand{\cH}{\mathcal H}
\newcommand{\G}{\mathbb G}

\newcommand{\cV}{\mathcal V}
\newcommand{\cC}{\mathcal C}
\newcommand{\mbX}{\mathbf{X}}
\newcommand{\bV}{\mathbb V}
\newcommand{\bW}{\mathbb W}
\newcommand{\m}{\mathrm m}
\newcommand{\cg}{\mathfrak g}

\newcommand{\bH}{\mathbb H}

\newcommand{\cG}{\mathcal G}
\newcommand{\res}{\mathbin{\llcorner}}
\newcommand{\Lie}{\mathrm{Lie}}
\newcommand{\beq}{\begin{equation}}
\newcommand{\eeq}{\end{equation}}
\newcommand{\cL}{\mathcal L}

\DeclareMathOperator{\Tan}{Tan}
\newcommand{\newsharedtheorem}[2]{%
    \IfFormatAtLeastTF{2026-06-01}{%
        \newtheorem{#1}[theorem]{#2}%
    }{%
        \newaliascnt{#1}{theorem}%
        \newtheorem{#1}[#1]{#2}%
        \aliascntresetthe{#1}%
    }%
}
\theoremstyle{plain}
\newtheorem{theorem}{Theorem}[section]
\newsharedtheorem{lemma}{Lemma}
\newsharedtheorem{proposition}{Proposition}
\newsharedtheorem{corollary}{Corollary}

\theoremstyle{definition}
\newsharedtheorem{definition}{Definition}

\theoremstyle{remark}
\newsharedtheorem{remark}{Remark}

\theoremstyle{definition}
\newsharedtheorem{example}{Example}
\newtheorem*{workedexample}{Example}
\newtheorem*{problem}{Problem}

\crefname{theorem}{Theorem}{Theorems}
\Crefname{theorem}{Theorem}{Theorems}
\crefname{lemma}{Lemma}{Lemmas}
\Crefname{lemma}{Lemma}{Lemmas}
\crefname{proposition}{Proposition}{Propositions}
\Crefname{proposition}{Proposition}{Propositions}
\crefname{corollary}{Corollary}{Corollaries}
\Crefname{corollary}{Corollary}{Corollaries}
\crefname{definition}{Definition}{Definitions}
\Crefname{definition}{Definition}{Definitions}
\crefname{remark}{Remark}{Remarks}
\Crefname{remark}{Remark}{Remarks}
\crefname{example}{Example}{Examples}
\Crefname{example}{Example}{Examples}
\newcommand{\proofstep}[1]{\par\medskip\noindent\textbf{Step~#1.}\enspace}

\title[$C^{m,\alpha}$-rectifiability in Carnot groups]
{$C^{m,\alpha}$-rectifiability in Carnot groups}

\author[K. O. Idu]{Kennedy Obinna Idu}
\address[Kennedy Obinna Idu]{University of Toronto, Department of Mathematics, Toronto, ON M5S 2E4, Canada.}
\email{idu@math.toronto.edu}
\author[F. P. Maiale]{Francesco Paolo Maiale}
\address[Francesco Paolo Maiale]{Gran Sasso Science Institute, Viale F. Crispi 7, 67100 L'Aquila, Italy.}
\email{francescopaolo.maiale@gssi.it}


\date{}

\begin{document}

\begin{abstract}
In an arbitrary Carnot group $\G$ of homogeneous dimension $Q$, we characterize rectifiability of sets of finite Hausdorff measure by $\mathbf C_\G^{m,\alpha}$-regular surfaces with complemented vertical tangent subgroups, for every $m\ge1$ and $\alpha\in[0,1]$. The criterion combines positive lower density with approximation in measure by tubes around such surfaces. At almost every point, one tube radius $\lambda r^{m+\alpha}$ suffices for $\alpha>0$, whereas the estimate must hold for every $\lambda r^m$, $\lambda>0$, for $\alpha=0$. Tubular approximation upgrades $\mathbf C^1_\G$-rectifiability to class $(m,\alpha)$ without an additional density assumption.

The proof combines ambient normal-form polynomial jets, a Remez-type norming inequality, and Whitney extension. In topological dimension at least two, quantitative volume approximation by noncharacteristic sublevel sets yields higher-order rectifiability of reduced boundaries. In particular, for weak solutions of the horizontal Poisson equation with $L^q_{\mathrm{loc}}$ forcing ($q>\max\{Q,2\}$ finite), almost every sublevel set has second-order rectifiable reduced boundary.
\end{abstract}

\maketitle

\section{Introduction} \label{sec:introduction}

Extending the classical theory of rectifiable sets beyond Euclidean spaces has been a central theme of geometric measure theory; see~\cites{AmbKir2000Rect,kirchheim1994rectifiable} and, for connections with harmonic analysis and singular integrals,~\cites{chousionis2019boundednes,DS91,DS93}. A rectifiable set is one that, although possibly very irregular, decomposes up to a null set into pieces modeled on regular surfaces. Federer's definition~\cite{Federer69}*{Section~3.2.14}, by almost-everywhere coverings with Lipschitz images, remains meaningful in metric spaces but does not generally capture intrinsic regular surfaces in Carnot groups. The equivalent Euclidean description by Lipschitz graphs is not directly available there; see~\cites{AmbKir2000Rect,idumama1}. A natural alternative, modeled on smooth submanifolds as non-critical level sets, is given by the \emph{regular surfaces} of Franchi, Serapioni, and Serra Cassano~\cites{FSSC01,FSSC6,FSSC11}, which are well suited to the low-codimensional regime. Different notions of rectifiability in Carnot groups do not coincide in general, because stratified geometry interacts nontrivially with projections, intrinsic graphs, and Whitney-type extension problems~\cites{FSSC4,pauls2004notion,AntLeD19,FranchiSerapioni2016IntrLip}. This motivates the following problem.

\begin{problem}[Rectifiability characterization]
Find necessary and sufficient geometric and \emph{quantitative} criteria for rectifiability of arbitrary order in metric spaces and, in particular, in Carnot groups.
\end{problem}

Quantitative criteria are scale-invariant estimates capturing approximation at small scales, often formulated through $\beta$-numbers~\cites{ninidu1,Chousionis2019}. In Euclidean spaces the problem goes back to Anzellotti and Serapioni~\cite{AnzSer}, and recent work gives quantitative criteria, for $\mathbf C^{1,\alpha}$-rectifiability in~\cite{ninidu1} and at every order in~\cite{San19}; see the survey~\cite{Mat23_Survey}. Complementary Euclidean criteria for $\mathbf C^{1,\alpha}$-rectifiability use square functions~\cite{Ghinassi2020} and averaged discrete curvatures~\cite{Kolasinski2017}. Approximate tangent characterizations are classical~\cite{mattila}*{Theorem~15.19}, and the sharp higher-order benchmark is~\cite{San19}*{Theorem~5.4}. In the present notation, it proves countable $\mathbf C^{m,\alpha}$-rectifiability, under $\sigma$-finiteness of the relevant Hausdorff measure, from approximation outside every tube of radius $\eps r^m$ and, when $\alpha>0$, outside one tube of radius $\lambda r^{m+\alpha}$, with no lower-density hypothesis.

In Heisenberg groups, first-order approximate-tangent characterizations were proved in~\cite{MSSC10}*{Theorems~3.14 and~3.15}. For $0<\alpha\le1$ and $n<k\le2n$,~Idu and Maiale characterize low-codimensional $\mathbf C_H^{1,\alpha}$-rectifiability by approximate tangent paraboloids under positive lower density~\cite{IM24}*{Theorem~1.1 and Proposition~1.2}. Their Remark~5.2 observes that the argument extends to codimension one in arbitrary Carnot groups, since its one exceptional lemma becomes immediate when the horizontal complement is one-dimensional. Complementary first-order results include Merlo's Marstrand--Mattila criterion for measures of codimension one~\cite{Merlo2023Marstrand} and the representation and density theorems of Antonelli and Merlo for measures with positive lower density, finite upper density, and a unique flat tangent up to multiples~\cites{AntonelliMerlo2022Representation,AntonelliMerlo2022Density}.

The main contribution of this paper is a characterization of intrinsic rectifiability of every order $m\ge1$ and every exponent $\alpha\in[0,1]$ by tubular approximation and positive lower density (\cref{thm:main,thm:main_reverse}). The result holds in arbitrary Carnot groups, in every dimension admitting vertical complemented tangent subgroups. Relative to the $\mathbf C_H^{1,\alpha}$ criterion of~\cite{IM24}, it treats both arbitrary orders and the full range of these tangent geometries. It provides a Carnot-group counterpart of Santilli's higher-order Euclidean criterion~\cite{San19}, with an additional positive lower-density hypothesis. That hypothesis is used only at first order: the higher-order step, \cref{thm:higher_order}, needs no density assumption.

In Carnot groups, surfaces with vertical tangent subgroups have two features absent in Euclidean spaces. First, the tangent objects are vertical homogeneous subgroups, which are graded but need not be Carnot (\cref{rmk:graded_not_carnot}), so sub-Riemannian Whitney and Lusin tools cannot in general be applied directly on a parameter subgroup $\bW$. Second, if a $\G$-regular surface is an intrinsic graph over $\bW$, the graph function need not be regular along the left-invariant fields of $\bW$: by the Baker--Campbell--Hausdorff formula, differentiating $w\mapsto w\cdot\phi(w)$ in a fixed $\bW$-direction may differentiate the defining function in higher-layer directions. The central higher-order step is to recover compatible ambient jets from approximation in measure, without differentiating the intrinsic graph parametrizations. On first-order graph patches, normal-form polynomials encode the approximating surfaces, and a norming inequality converts control of polynomial values on subsets of positive density into control of all ambient derivatives through order $m$. Applied to differences of the local polynomials, this yields the Whitney compatibility estimates needed to construct rectifying surfaces as non-critical level sets of ambient $\mathbf C_\G^{m,\alpha}$ functions.

At $\alpha=1$ the conclusion is class $(m,1)$, not $(m+1,0)$. In $\R^n$ the two notions of rectifiability coincide~\cite{AnzSer}*{Proposition~3.2}, since a $C^{m,1}$ graph agrees with a $C^{m+1}$ graph outside a set of small measure. In positive codimension in Carnot groups, this Lusin-type step would have to be run on the top derivatives of the extension restricted to the surface, which is Haar-null, so Pansu differentiability almost everywhere gives no information there, and the graph maps of intrinsic parametrizations need not be regular along the parameter subgroup; it would require an extension theorem for adapted trace jets beyond the Whitney tools used here. In Euclidean space, Brena, De Lellis, and Franceschini~\cite{BrenaDeLellisFranceschini2025} prove $\mathbf C^\infty$-rectifiability of the supports of stationary integral varifolds of dimension at least two. Related results concern second-order rectifiability of Euclidean varifolds~\cites{MR3023856,MR4243013} and rectifiability of singular sets of area-minimizing currents~\cite{MR4888467}; higher-order criteria of this type may also be useful for sub-Riemannian variational problems.

\subsection{Main results} \label{subsec:main_result}

Let $N:=\dim\G$, let $Q$ be the homogeneous dimension, let $1\le k\le N$, and set $k_{\m}:=Q-(N-k)$. The tangent subgroups of $k$-dimensional $\G$-regular surfaces (\cref{def:regular_surface}) are vertical and admit a homogeneous complement, that is, they belong to the set $\cV(\G,k)$ of~\cref{subsec:grassmannian}. Such surfaces exist if and only if $\cV(\G,k)\neq\emptyset$, and otherwise all statements below are empty (\cref{rmk:beyond_low_codim}). In the low-codimensional range $\max\{1,N-\upsilon\}\le k\le N$ of~\cref{def:upsilon}, every complemented $k$-dimensional homogeneous subgroup is vertical. The threshold describes the type of subgroup whenever one exists. The family $\cV(\G,k)$ can also be nonempty below this range. We write $\Theta_*^s(A,p)$, $\Theta^{*s}(A,p)$, and $\Theta^s(A,p)$ for the lower, upper, and exact $s$-dimensional densities of $\cH^s\res A$ at $p$ (see~\cref{sec:preliminaries} for more details).

The tubular approximation is measured in normalized $\cH^{k_{\m}}$-measure outside a tube around an ambient $\G$-regular surface, of radius of order $r^{m+\alpha}$ when $\alpha>0$ and, when $\alpha=0$, of radius $\lambda r^m$ for every $\lambda>0$ (a measure-theoretic little-$o(r^m)$ condition). %No intrinsic-graph representation of the approximating surface is required.

\begin{definition} \label{def:special_class}
Let $\alpha\in[0,1]$, $m\in\N$, and $1\le k\le N$. A set $E\subset\G$ belongs to the class $\mathfrak{S}_{m,\alpha}(k)$ if $E$ is $\cH^{k_{\m}}$-measurable, $\cH^{k_{\m}}(E)<\infty$, and for $\cH^{k_{\m}}$-almost every $p\in E$ there is a $k$-dimensional $\G$-regular surface $S=S_p$ of class $(m,\alpha)$, in the sense of~\cref{def:regular_surface}, with $p\in S$ and the following property:
\begin{enumerate}[label=(\alph*)]
    \item if $\alpha>0$, there is a constant $\lambda=\lambda(p)\ge0$ such that
    \begin{equation}\label{eq:main_assumption}
    \lim_{r \downarrow 0} \frac{1}{r^{k_{\m}}} \cH^{k_{\m}} \Bigl( E \cap B_r(p) \cap \bigl\{ q \in \G \ : \ d(q,S) > \lambda\, r^{m+\alpha} \bigr\} \Bigr) = 0;
    \end{equation}
    \item if $\alpha=0$, then for every $\lambda>0$,
    \begin{equation}\label{eq:main_assumption_alpha_zero}
    \lim_{r \downarrow 0} \frac{1}{r^{k_{\m}}} \cH^{k_{\m}} \Bigl( E \cap B_r(p) \cap \bigl\{ q \in \G \ : \ d(q,S) > \lambda\, r^m \bigr\} \Bigr) = 0.
    \end{equation}
\end{enumerate}
\end{definition}

\begin{theorem} \label{thm:main}
Let $\alpha\in[0,1]$, $m\in\N$, and $1\le k\le N$. If $E\in\mathfrak{S}_{m,\alpha}(k)$ and $\Theta_*^{k_{\m}}(E,p)>0$ for $\cH^{k_{\m}}$-almost every $p\in E$, then $E$ is $(k,\G)$-rectifiable of class $(m,\alpha)$ with respect to $\cH^{k_{\m}}$, in the sense of~\cref{def:rec}.
\end{theorem}

The higher-order part of~\cref{thm:main} requires no density assumption once the first-order rectifiability is settled.

\begin{theorem} \label{thm:higher_order}
Let $\alpha\in[0,1]$, $m\in\N$, and $1\le k\le N$. If $E\in\mathfrak{S}_{m,\alpha}(k)$ is $(k,\G)$-rectifiable of class $(1,0)$ with respect to $\cH^{k_{\m}}$, then $E$ is $(k,\G)$-rectifiable of class $(m,\alpha)$ with respect to $\cH^{k_{\m}}$.
\end{theorem}

\begin{remark}\label{rmk:density_order_one}
By~\cref{thm:higher_order}, positive lower density enters~\cref{thm:main}
only through~\cref{thm:c1_rec}, where it turns approximate tangency in
measure into uniform pairwise flatness on compact pieces. In Euclidean space
no density assumption is needed, at order one~\cite{mattila}*{Theorem~15.19}
or higher~\cite{San19}*{Theorem~5.4}; at order one, the set is covered by
Lipschitz graphs, which are $C^1$-rectifiable by Rademacher's theorem. An
intrinsic Rademacher theorem holds in Heisenberg groups~\cite{Vittone2022},
but some Carnot groups carry intrinsic Lipschitz graphs over vertical
subgroups that are nowhere intrinsically differentiable~\cite{JNGV2021Nowhere}.
Even in Heisenberg groups, whether the density assumption
of~\cite{MSSC10}*{Theorem~3.15} is needed was recorded as open
in~\cite{Mat23_Survey}*{Theorem~8.7 and the remark after it} and recently resolved in~\cite{Idu2025Density}. The authors are investigating the corresponding density question in arbitrary Carnot groups.
\end{remark}

The converse requires no hypothesis beyond rectifiability and finite measure.

\begin{theorem} \label{thm:main_reverse}
Let $\alpha\in[0,1]$, $m\in\mathbb N$, and $1\le k\le N$. If $E\subset\G$ is $(k,\G)$-rectifiable of class $(m,\alpha)$ and $\cH^{k_{\m}}(E)<\infty$, then $E\in\mathfrak{S}_{m,\alpha}(k)$ and
\begin{equation}\label{eq:density_assumption}
    \Theta_*^{k_{\m}}(E,p)>0\qquad\text{for $\cH^{k_{\m}}$-almost every $p\in E$.}
\end{equation}
\end{theorem}

\begin{proof}
If $k=N$, take $S_p=\G$. Let $k<N$ and $s:=k_{\m}$. By~\cref{def:rec} there are
$k$-dimensional $\G$-regular surfaces $S_i$ of class $(m,\alpha)$ with
$\cH^s(E\setminus\bigcup_iS_i)=0$, and the sets
$A_i:=E\cap S_i\setminus\bigcup_{j<i}S_j$ are pairwise disjoint and cover $E$
up to a null set. Each $S_i$ is Borel, being a countable union of sets of the
form $f^{-1}(0)\cap\Omega$, and $\mu_i:=\cH^s\res(E\setminus S_i)$ is a finite
Radon measure with $\mu_i(S_i)=0$. The density estimate
of~\cite{AGV2022}*{Theorem~2.1(ii)}, applied as in the proof of its
Corollary~3.6, therefore gives zero upper Federer $s$-density of $\mu_i$
at $\cH^s$-almost every point of $S_i$; see also~\cite{Federer69}*{2.10.19(4)}. Since centered closed balls are
among the sets defining that density and have diameter at most $2r$,
we obtain
\[
    \lim_{r\downarrow0}r^{-s}\cH^s\bigl((E\setminus S_i)\cap B_r(p)\bigr)=0
    \qquad\text{for }\cH^s\text{-a.e.\ }p\in A_i.
\]
At such $p$ take $S=S_i$, and $\lambda=0$ if $\alpha>0$. The exceptional sets
in~\eqref{eq:main_assumption} and~\eqref{eq:main_assumption_alpha_zero} are
then contained in $(E\setminus S_i)\cap B_r(p)$, so
$E\in\mathfrak{S}_{m,\alpha}(k)$.

For the density assertion, if $k=N$ the claim follows from the Lebesgue
density theorem for Haar measure. Let $k<N$. Refining the cover to
relatively compact surface patches, we may assume that
$\nu_i:=\cH^s\res S_i$ is finite. By
\cite{AGV2022}*{Propositions~2.3(ii)--(iii) and~2.17} (stated there for
closed balls; the open-ball version follows since $B_r\subset\overline B_r$
and $\nu_i(B_r(p))\ge\nu_i(\overline B_{r/2}(p))$, which is comparable to
$\nu_i(\overline B_r(p))$ by local Ahlfors regularity),
\[
 \frac{\nu_i((S_i\setminus E)\cap B_r(p))}{\nu_i(B_r(p))}\longrightarrow0
 \qquad\text{for }\cH^s\text{-a.e.\ }p\in E\cap S_i.
\]
At such $p$ the local lower Ahlfors bound gives
$\cH^s(E\cap B_r(p))\ge\tfrac12\nu_i(B_r(p))\ge c_pr^s$ for all small $r$.
Alternatively, the existence and positivity of the density follow
from~\cite{AGV2022}*{Corollary~3.6}, applied to $E\cap S_i$.
\end{proof}

Together with~\cref{thm:main}, this shows that an $\cH^{k_{\m}}$-measurable set
$E$ of finite $\cH^{k_{\m}}$-measure is rectifiable of class $(m,\alpha)$
if and only if $E\in\mathfrak S_{m,\alpha}(k)$ and positive lower
density~\eqref{eq:density_assumption} holds almost everywhere.

In \cref{sec:finite_perimeter} we apply these results to sets of locally
finite perimeter in an arbitrary Carnot group of dimension at least two: the
reduced boundary is rectifiable of class $(m,\alpha)$ if and only if it
satisfies the tubular condition (\cref{cor:reduced_boundary}), since positive
lower density is automatic (\cref{lem:perimeter_facts}),
and this holds whenever the set is approximated in volume, to order
$o(r^{Q-1+m+\alpha})$, by sublevel sets of $\mathbf C^{m,\alpha}_\G$ functions
with nonvanishing horizontal gradient, for instance polynomials
(\cref{lem:bulk_blowups}). For perimeter quasiminimizers, the topological
boundary of the open representative differs from the reduced boundary by a
null set, so these conclusions apply to it as well (\cref{rmk:bulk_examples}).

As a concrete PDE application, for weak solutions of the horizontal Poisson
equation (for instance, Dirichlet minimizers) sufficiently integrable forcing yields
the quadratic bulk rate on the reduced boundaries of almost every sublevel set,
giving second-order rectifiability (\cref{prop:dirichlet_levels}). Finally, from
order equal to the step of $\G$ on, the covering regular surfaces are also
Euclidean $C^1$ submanifolds (\cref{sec:euclidean_regularity}).

\subsection*{Organization of the paper}

\Cref{sec:preliminaries,sec:rect_tech} recall the Carnot-group framework,
regular surfaces, and the approximate-tangent criterion. \Cref{sec:whitney_jets}
develops the Whitney-extension, blow-up, normal-form, and jet-selection tools,
which \cref{sec:higher-order-rectifiability} assembles into the proofs
of~\cref{thm:c1_rec,thm:higher_order,thm:main}.

%%%%%%%%%%%%%%%%%%%%%%%%%%%%%%%%%%%%%%%%%%%%%%%%%%%%%%%%%%%%%%%%%%%%%%%%%%%%%
\section{Preliminaries}\label{sec:preliminaries}
%%%%%%%%%%%%%%%%%%%%%%%%%%%%%%%%%%%%%%%%%%%%%%%%%%%%%%%%%%%%%%%%%%%%%%%%%%%%%

\subsection{Graded and Carnot groups} \label{subsec:graded_Lie}

A simply connected nilpotent Lie group $\G$ is a \emph{graded group} if its Lie algebra admits a decomposition
\begin{equation}\label{eq:graded_algebra}
\cg = \bigoplus_{i=1}^\iota \cg_i, \qquad [\cg_i, \cg_j] \subseteq \cg_{i+j},
\end{equation}
with $\cg_j=\{0\}$ for $j>\iota$ and $\cg_\iota\ne\{0\}$. It is a \emph{Carnot} (or \emph{stratified}) group of step $\iota$ if, in addition, it satisfies
\begin{equation}\label{eq:stratified_algebra}
[\cg_1, \cg_j] = \cg_{j+1}, \qquad j=1,\ldots,\iota-1.
\end{equation}
Set $h_i:=\dim\cg_1+\cdots+\dim\cg_i$, $h_0:=0$, $\ell:=h_1$, $N:=h_\iota=\dim\cg$, and $Q:=\sum_{i=1}^\iota i\dim\cg_i$ to denote the homogeneous dimension. 

We fix a basis $X_1,\ldots,X_N$ adapted to~\eqref{eq:graded_algebra}, i.e. $X_{h_{i-1}+1},\ldots,X_{h_i}$ is a basis of $\cg_i$, and give the $l$-th coordinate the degree $d_l:=j$ when $h_{j-1}<l\le h_j$. The exponential map is a global diffeomorphism~\cite{BLU2007stratified}*{Theorem~1.3.28}, and we identify $\G$ with $\R^N$ through the exponential coordinates $x\mapsto\exp(x_1X_1+\cdots+x_NX_N)$. Then $e=0$, $x^{-1}=-x$, and $\R^N=\G^1\oplus\cdots\oplus\G^\iota$, where the layer $\G^i=\exp(\cg_i)$ need not be a subgroup. We write $|\cdot|$ for the Euclidean norm in which the adapted basis is orthonormal. The \emph{dilations} $\delta_\lambda(x)=(\lambda^{d_1}x_1,\ldots,\lambda^{d_N}x_N)$, $\lambda>0$, are group automorphisms.

\begin{remark}\label{rmk:BCH}
By the Baker--Campbell--Hausdorff (BCH) formula,
\begin{equation}\label{eq:BCH}
x \cdot y = x + y + \tfrac{1}{2}[x,y] + \sum_{j=3}^{\iota} c_j(x,y),
\end{equation}
where $c_j$ has bracket length $j$; see~\cite{BLU2007stratified}*{Theorem~1.3.28}. Since dilations are automorphisms, we have $(x\cdot y)_l=x_l+y_l+P_l(x,y)$, where $P_l$ is a polynomial, homogeneous of degree $d_l$ in $(x,y)$, depending only on coordinates of degree strictly less than $d_l$ (in particular, $P_l\equiv0$ when $d_l=1$) and satisfying $P_l(x,0)=P_l(0,y)=0$.
\end{remark}

We fix a \emph{homogeneous distance} $d$, that is, a left-invariant distance inducing the Euclidean topology and satisfying $d(\delta_\lambda x,\delta_\lambda y)=\lambda d(x,y)$, and set $\|x\|:=d(x,e)$; then
\[
\|\delta_\lambda x\|=\lambda\|x\|, \qquad \|x^{-1}\|=\|x\|, \qquad \|xy\|\le\|x\|+\|y\|. 
\]
Any two homogeneous distances are bi-Lipschitz equivalent. Their norms have a positive continuous ratio on the compact set $\{x:\sum_l|x_l|^{1/d_l}=1\}$, and homogeneity extends these bounds to $\G$. Let $B_r(x)$ and $\overline B_r(x)$ be the open and closed balls, $d(x,S):=\inf_{y\in S}d(x,y)$, and $\cH^s$ the $s$-dimensional Hausdorff measure with respect to $d$; see~\cite{Federer69}*{Section~2.10.2}. For a Radon measure $\mu$ and $s\ge0$ we set
\[
\Theta_*^s(\mu,x):=\liminf_{r\downarrow0}\frac{\mu(B_r(x))}{r^s},\qquad
\Theta^{*s}(\mu,x):=\limsup_{r\downarrow0}\frac{\mu(B_r(x))}{r^s},
\]
for the lower/upper densities, and write $\Theta^s(\mu,x)$ for their common value when they agree. We abbreviate $\Theta_*^s(A,x):=\Theta_*^s(\cH^s\res A,x)$, and similarly for $\Theta^{*s}$ and $\Theta^s$. The Lebesgue measure $\cL^N$ is a Haar measure of $\G$ (\cite{BLU2007stratified}*{1.3.21}), and $\cL^N(B_r(p))=r^Q\cL^N(B_1(e))$.

For a multi-index $J\in\N_0^N$, let $|J|:=j_1+\cdots+j_N$ and $d(J):=\sum_{i=1}^N d_ij_i$ be its length and degree respectively. 

\begin{definition}
A polynomial has \emph{homogeneous degree at most $m$} if it is a linear combination of monomials $x^J$ with $d(J)\le m$, and is \emph{homogeneous of degree $m$} if only monomials with $d(J)=m$ occur.
\end{definition}

Identifying each $X_i$ with the corresponding left-invariant vector field on $\G$, by~\cite{FS82}*{Section~1C} there are homogeneous polynomials $\mathcal Q_{ik}$ of degree $d_k-d_i$ such that
\[
X_i = \partial_{x_i}+\sum_{d_k>d_i}\mathcal Q_{ik}(x)\,\partial_{x_k}.
\]
We set $\mbX^J:=X_1^{j_1}\cdots X_N^{j_N}$, a left-invariant differential operator, homogeneous of degree $d(J)$ with respect to dilations.

%%%%%%%%%%%%%%%%%%%%%%%%%%%%%%%%%%%%%%%%%%%%
\subsection{Projections and Grassmannian}\label{subsec:grassmannian}
%%%%%%%%%%%%%%%%%%%%%%%%%%%%%%%%%%%%%%%%%%%%

Recall that a \emph{homogeneous subgroup} of $\G$ is a closed Lie subgroup $\mathbb L$ with the property that $\delta_\lambda\mathbb L=\mathbb L$ for all $\lambda>0$.

\begin{remark}\label{rmk:graded_not_carnot}
Homogeneous subgroups are linear subspaces in exponential coordinates and
inherit the grading. Indeed, $\mathbb L$ is closed, hence an embedded
submanifold, so the paths $t\mapsto\delta_t p$, extended by $e$ at $t=0$,
are continuous in $\mathbb L$ and connect every $p\in\mathbb L$ to $e$. Thus
$\mathbb L$ is connected and $\mathbb L=\exp(\mathfrak h)$
with $\mathfrak h:=\Lie(\mathbb L)$. The space $\mathfrak h$ is invariant under
$D\delta_2(e)$ and hence under its spectral projections onto the layers
$\cg_i$. Thus $\mathfrak h=\bigoplus_{i=1}^\iota(\mathfrak h\cap\cg_i)$ and
$\mathbb L=\mathbb L^1\oplus\cdots\oplus\mathbb L^\iota$ with $\mathbb L^i\subseteq\G^i$.
Every homogeneous subgroup is therefore a graded group, with the restricted
dilations and BCH product~\eqref{eq:BCH}. However, it \emph{need not be Carnot}:
$\mathbb L^1$ may fail to generate the higher layers of $\Lie(\mathbb L)$. This is
one of the main issues that we have to overcome in this paper.
\end{remark}

\begin{remark}\label{rmk:haar_measure}
If $\mathbb L\le\G$ is a homogeneous subgroup with homogeneous dimension $Q_{\mathbb L}:=\sum_{i=1}^\iota i\dim\mathbb L^i$, its Haar measure $\mu_{\mathbb L}$ satisfies $\mu_{\mathbb L}(\delta_\lambda A)=\lambda^{Q_{\mathbb L}}\mu_{\mathbb L}(A)$. Ball scaling makes $\mu_{\mathbb L}$ Ahlfors $Q_{\mathbb L}$-regular, so $\cH^{Q_{\mathbb L}}\res\mathbb L$ is nonzero and locally finite. It is left invariant, hence a constant multiple of $\mu_{\mathbb L}$ by uniqueness of Haar measure.
\end{remark}

A homogeneous subgroup $\mathbb L$ is \emph{horizontal} if $\mathbb L=\mathbb L^1$, and \emph{vertical} if $\mathbb L^j=\G^j$ for all $j\ge2$; vertical subgroups are normal. 

\begin{definition}[Splitting]
Two homogeneous subgroups $\bW,\bV$ are \emph{complementary} if $\bW\cap\bV=\{0\}$ and $\G=\bW\cdot\bV$; if $\bW$ is normal, we write $\G=\bW\rtimes\bV$. 
\end{definition}

For $0\le k\le N$, let $\cG(\G,k)$ be the set of $k$-dimensional homogeneous subgroups admitting a complementary homogeneous subgroup, let $\cV(\G,k)\subseteq\cG(\G,k)$ denote its vertical elements, and set $\cG(\G):=\bigcup_{k=0}^N\cG(\G,k)$ for the \emph{intrinsic Grassmannian}. The following standard projection facts are discussed in~\cites{MSSC10,FranchiSerapioni2016IntrLip}; see~\cite{MSSC10}*{Proposition~2.15} in Heisenberg groups and, in homogeneous groups,~\cite{idumama1}*{Theorem~2.13} for the splitting $\G=\bV^\perp\rtimes\bV$ of a horizontal subgroup $\bV$ and its orthogonal complement $\bV^\perp$.

\begin{theorem} \label{thm:semidirect_fact}
	If $\G = \bW \rtimes \bV$ with $\bV$ horizontal and $\bW$ vertical, each $x \in \G$ has unique components $x=\pi_\bW(x)\cdot\pi_\bV(x)$ with $\pi_\bW(x)\in\bW$ and $\pi_\bV(x)\in\bV$, and
	\begin{equation} \label{eq:projection_algebra}
		\begin{aligned}
			\pi_\bW(x^{-1}) &= \pi_\bV(x)^{-1} \pi_\bW(x)^{-1} \pi_\bV(x), &\quad& \pi_\bV(x^{-1}) = \pi_\bV(x)^{-1}, \\
			\pi_\bW(\delta_\lambda x) &= \delta_\lambda \pi_\bW(x), &\quad& \pi_\bV(\delta_\lambda x) = \delta_\lambda \pi_\bV(x), \\
			\pi_\bW(xy) &= \pi_\bW(x)\pi_\bV(x)\pi_\bW(y)\pi_\bV(x)^{-1}, &\quad& \pi_\bV(xy) = \pi_\bV(x)\pi_\bV(y).
		\end{aligned}
	\end{equation}
	Moreover, the constant
	$C_{\bW,\bV}:=\max\bigl\{1,\sup_{\|x\|=1}\bigl(\|\pi_\bV(x)\|+\|\pi_\bW(x)\|\bigr)\bigr\}$
	is finite and, for all $x\in\G$,
	\begin{equation}\label{eq:piVcontin}
		\|\pi_\bV(x)\| + \|\pi_{\bW}(x)\| \le C_{\bW,\bV}\|x\|,
	\end{equation}
	\begin{equation}\label{eq:projection_distance}
		\begin{aligned}
			& C_{\bW,\bV}^{-1}\|\pi_\bV(x)\| \le d(x, \bW) \le \|\pi_\bV(x)\|, \\
			& C_{\bW,\bV}^{-1}\|\pi_\bV(x)^{-1}\pi_{\bW}(x)\pi_\bV(x)\| \le d(x, \bV) \le \|\pi_\bV(x)^{-1}\pi_{\bW}(x)\pi_\bV(x)\|.
		\end{aligned}
	\end{equation}
\end{theorem}

\begin{proof}
The identities follow from the uniqueness of the factors.
By~\eqref{eq:BCH}, first-layer coordinates add under products, so
$\G=\bW\cdot\bV$ gives $\cg_1\subseteq\Lie(\bW)+\Lie(\bV)$, while the higher layers
lie in $\Lie(\bW)$ because $\bW$ is vertical. Since $\exp$ is injective and
$\bW\cap\bV=\{0\}$, the sum is direct, i.e. $\cg=\Lie(\bW)\oplus\Lie(\bV)$. As
$\bW$ is normal, the
multiplication map has invertible differential everywhere and is bijective, so its
inverse, and hence both projections, are smooth. Homogeneity and compactness
of the unit sphere give~\eqref{eq:piVcontin}. For $w\in\bW$, apply that bound
to $w^{-1}x$ to get the lower bound for $d(x,\bW)$, while the upper bound follows
by choosing $w=\pi_\bW(x)$. For $v\in\bV$,
$\pi_\bW(x^{-1}v)=\pi_\bV(x)^{-1}\pi_\bW(x)^{-1}\pi_\bV(x)$ gives
the other lower bound, while choosing $v=\pi_\bV(x)$ gives its upper bound.
\end{proof}

By~\eqref{eq:projection_algebra}, $\pi_\bV:\G\to\bV$ is a homomorphism commuting with dilations, that is, an h-homomorphism in the sense of~\cref{def:h-homomorphism}. If $\bV$ is $c$-dimensional and a basis of $\bV$ is fixed, we write $\Pi_{\bV}:=\operatorname{coord}_{\bV}\circ\, \pi_{\bV}:\G\to\R^c$ for the coordinate vector of the group element $\pi_\bV(x)$. This distinction between projections matters in~\cref{sec:whitney_jets}.

\begin{definition}[Threshold] \label{def:upsilon}
	Let 
    \[
    t_\G := \min\bigl\{s\in\{0,\ldots,N\} \ : \ \cG(\G,k)=\cV(\G,k) \text{ for every } s \le k \le N \bigr\}.
    \]
    The \emph{low-codimension threshold} of $\G$ is $\upsilon := N-t_\G$; thus, all $k$-dimensional elements of the intrinsic Grassmannian are vertical whenever $N-\upsilon \le k \le N$.
\end{definition}

The value of $\upsilon$ depends on the group, as the examples below show. They rely on the following characterization of the vertical subgroups that admit a complement.

\begin{proposition}\label{prop:vertical_splitting}
If $S_1\subseteq\cg_1$ is a linear subspace and $\mathfrak w:=S_1\oplus\cg_2\oplus\cdots\oplus\cg_\iota$, then the vertical subgroup $\bW=\exp(\mathfrak w)$ admits a complementary homogeneous subgroup if and only if $\cg_1=S_1\oplus\mathfrak h$ for some abelian subspace $\mathfrak h$, and then $\G=\bW\rtimes\exp(\mathfrak h)$. Moreover, every homogeneous complement of $\bW$ is of this form, i.e. horizontal and abelian.
\end{proposition}

\begin{proof}
Let $\bV=\exp(\mathfrak h)$ be a homogeneous complement of $\bW$. By~\cref{rmk:graded_not_carnot}, $\mathfrak h$ is graded; since $\mathfrak h\cap\mathfrak w=\{0\}$ and $\mathfrak w$ contains every layer $\cg_j$ with $j\ge2$, we get $\mathfrak h\subseteq\cg_1$. The multiplication map $\bW\times\bV\to\G$ is a smooth surjection, so $\dim\bW+\dim\bV\ge N$ by Sard's theorem. Hence $\cg=\mathfrak w\oplus\mathfrak h$ and $\cg_1=S_1\oplus\mathfrak h$. Moreover, $\mathfrak h$ is a subalgebra, so $[\mathfrak h,\mathfrak h]\subseteq\mathfrak h\cap\cg_2=\{0\}$.

Conversely, let $\cg_1=S_1\oplus\mathfrak h$ with $\mathfrak h$ abelian. Then $\mathfrak h$ is a subalgebra, $\bV:=\exp(\mathfrak h)$ is a horizontal homogeneous subgroup, and $\bW$ is normal, so $\bW\bV$ is a subgroup. Since $\cg=\mathfrak w\oplus\mathfrak h$, the inverse function theorem shows that $\bW\bV$ contains a neighborhood of $e$; an open subgroup is also closed, so $\bW\bV=\G$ by connectedness. Finally, $\bW\cap\bV=\{0\}$ because $\exp$ is injective and $\mathfrak w\cap\mathfrak h=\{0\}$.
\end{proof}

\begin{example}
We now provide a few concrete examples:
\begin{itemize}
    \item In $\bH^n$, as proved in~\cite{MSSC10}*{Proposition~2.17}, one has $t_{\bH^n}=n+1$, so $\upsilon=n$ and the low-codimensional range is $n+1\le k\le2n+1$.
    
    \item Let $V=\operatorname{span}\{X_1,\ldots,X_4\}$ and
    $Z=\operatorname{span}\{Z_1,Z_2\}$. The six-dimensional H-type algebra
    $\cg=V\oplus Z$ has nonzero brackets
    \[
    [X_1,X_2]=[X_3,X_4]=Z_1,\qquad
    [X_1,X_3]=-[X_2,X_4]=Z_2,
    \]
    see, e.g.,~\cite{Balashchenko2011}*{(4.2)}. The bracket is
    $\omega_1Z_1+\omega_2Z_2$, where
    \[
    \omega_1=X_1^*\wedge X_2^*+X_3^*\wedge X_4^*,
    \qquad
    \omega_2=X_1^*\wedge X_3^*-X_2^*\wedge X_4^*.
    \]
    For $(a,b)\ne(0,0)$ the form $a\omega_1+b\omega_2$ has Pfaffian
    $a^2+b^2\ne0$, so no three-dimensional subspace of $V$ is isotropic
    for it. Thus $[S,S]=Z$ whenever $\dim S=3$. A nonvertical
    homogeneous subalgebra $S\oplus Z'$ has $\dim Z'\le1$, hence
    $\dim S\le2$ and total dimension $\le 3$.
    Therefore $\cG(\G,k)=\cV(\G,k)$ for $4\le k\le6$.

    Sharpness follows from the complementary subalgebras
    \[
    \mathfrak w_0=\operatorname{span}\{X_1,X_2,Z_1\},
    \qquad
    \mathfrak h_0=\operatorname{span}\{X_3,X_1+X_4,Z_1-Z_2\}.
    \]
    Indeed $[X_3,X_1+X_4]=Z_1-Z_2$, their horizontal parts are
    complementary in $V$, and their central parts are complementary in
    $Z$. The step-two BCH formula first determines the horizontal factors
    and then the central factors uniquely. Thus
    $\exp(\mathfrak w_0)\exp(\mathfrak h_0)=\G$ with trivial intersection,
    so $t_\G=4$ and $\upsilon=2$.

    \item In the free step-two group on $r\ge3$ generators, two independent
    horizontal vectors never commute. By~\cref{prop:vertical_splitting},
    vertical subgroups are therefore complemented only in dimensions
    $N-1$ and $N$. The subalgebra
    \[
    \mathfrak w=
    \operatorname{span}\bigl\{X_2,\ldots,X_r\}\oplus
    \operatorname{span}\{[X_i,X_j] \ : \ 1\le i<j\le r,\ (i,j)\ne(1,2)\bigr\}
    \]
    defines a nonvertical subgroup $\exp(\mathfrak w)$ of dimension $N-2$,
    complemented by $\exp\operatorname{span}\{X_1,[X_1,X_2]\}$, again by
    the step-two BCH formula.

    \noindent The one-dimensional complement of a codimension-one complemented
    subgroup lies in $\cg_1$ or $\cg_2$. In the second case the subgroup
    would contain $\cg_1$ and hence $[\cg_1,\cg_1]=\cg_2$, which is
    impossible. Thus the complement is horizontal and the subgroup is
    vertical, proving $\upsilon=1$.
\end{itemize}
\end{example}

\begin{remark}\label{rmk:beyond_low_codim}
The results of this paper use only that the relevant $k$-dimensional
subgroups are vertical and complemented, and they are stated for every
$1\le k\le N$. If $\bW\in\cV(\G,k)$ has the horizontal complement $\bV$,
then $\Pi_\bV$ is a smooth h-homomorphism onto $\R^{N-k}$ with kernel $\bW$,
so $\bW$ is a $k$-dimensional $\G$-regular surface of every class
(\cref{def:regular_surface}). Conversely, the tangent subgroups of such
surfaces lie in $\cV(\G,k)$, as shown after~\cref{def:tangent_group}. Hence the
statements are nonempty when $\cV(\G,k)\neq\emptyset$. The
threshold $\upsilon$ marks the range in which every element of $\cG(\G,k)$
is vertical, but $\cV(\G,k)$ can be nonempty below it. For instance, let
$\G=\bH^1\times\R$, with $\cg_1=\operatorname{span}\{X,Y,U\}$, $\cg_2=\operatorname{span}\{T\}$,
and $[X,Y]=T$ the only nonzero bracket. The horizontal subgroup
$\exp(\operatorname{span}\{X,U\})$ is complemented by the normal subgroup
$\exp(\operatorname{span}\{Y,T\})$, while every three-dimensional
homogeneous subgroup contains $T$, since $\cg_1$ is not a subalgebra; hence
$t_\G=3$ and $\upsilon=1$. Nevertheless,
$\exp(\operatorname{span}\{X,T\})\in\cV(\G,2)$, with the abelian complement
$\exp(\operatorname{span}\{Y,U\})$. Not every two-dimensional vertical subgroup
is complemented there: by~\cref{prop:vertical_splitting},
$\exp(\operatorname{span}\{U,T\})$ is not, because every two-dimensional
abelian subspace of $\cg_1$ contains $U$.
\end{remark}

\begin{remark}\label{rmk:grassmannian_metric}
A split vertical subgroup may have many horizontal complements, and the Euclidean orthogonal complement of its first layer need not be abelian. We metrize $\cV(\G,k)$ through first layers. If $P_{\bW^1}$ is the orthogonal projection of $\G^1$ onto $\bW^1$, we set
\begin{equation}\label{eq:grass_top}
	\rho\bigl(\bW_1,\bW_2\bigr):=\|P_{\bW_1^1}-P_{\bW_2^1}\|_{\operatorname{op}}, \qquad \bW_1,\bW_2\in\cV(\G,k).
\end{equation}
In the range $N-\upsilon\le k\le N$ this is also the topology used on $\cG(\G,k)=\cV(\G,k)$.
\end{remark}

By~\cref{rmk:haar_measure}, every $\bW\in\cV(\G,k)$ has Hausdorff dimension
\[
Q_\bW=\dim\bW^1+\sum_{j\ge2}j\dim\cg_j=Q-(N-k).
\]
For the vertical dimensions considered below, we therefore set
\begin{equation}\label{eq:metric_dimension}
	k_{\m} := Q -(N-k), \qquad 1\le k\le N.
\end{equation}

%%%%%%%%%%%%%%%%%%%%%%%%%%%%%%%%%%%%%%%%%%%%%%%%%%%%%%%%%%%%%%%%%%%%%%%%%%%%%%%%
\subsection{Intrinsically regular functions} \label{subsec:polynomials}
%%%%%%%%%%%%%%%%%%%%%%%%%%%%%%%%%%%%%%%%%%%%%%%%%%%%%%%%%%%%%%%%%%%%%%%%%%%%%%%%

From this subsection onward, the ambient group $\G$ is Carnot.
Its homogeneous subgroups may still be only graded.

\begin{definition} \label{def:alpha_regular}
	Let $\Omega\subseteq \G$ be open, $m \in \mathbb{N}_0$, and $\alpha \in [0,1]$. A function $f:\Omega\to\R$ is \emph{intrinsically $(m,\alpha)$-regular}, written $f\in \mathbf{C}_\G^{m,\alpha}(\Omega)$, if:
	\begin{itemize}
		\item for $m=0$, $f$ is continuous when $\alpha=0$, and locally $\alpha$-H\"older with respect to $d$ when $0<\alpha\le1$;
		\item for $m\ge1$, $\mbX^J f$ exists and is continuous on $\Omega$ for every $J$ with $d(J) \le m$ and, when $0<\alpha\le 1$, for every compact $K\Subset \Omega$ and every $J$ with $d(J)=m$,
		\[
		\sup_{x \neq y \in K} \frac{\bigl| \mbX^J f(x)- \mbX^J f(y) \bigr|}{d(x,y)^\alpha} < \infty.
		\]
	\end{itemize}
	In addition, we write $\mathbf C_\G^m(\Omega):=\mathbf C_\G^{m,0}(\Omega)$ and call its elements intrinsically $m$-regular functions. For $F=(F_1,\ldots,F_c)\in[\mathbf C_\G^1(\Omega)]^c$, we write $\nabla_HF(p):=(X_jF_i(p))_{i\le c,\,j\le\ell}$ for the horizontal Jacobian.
\end{definition}

A \emph{polynomial of homogeneous degree at most $m$ centered at $x_0$} is an expression of the form $\sum_{d(J)\le m}a_J(x_0)(x_0^{-1}x)^J/J!$. For $f \in \mathbf{C}_\G^m(\Omega)$ and $x_0 \in \Omega$, the \emph{Taylor polynomial} $P_m[f](x_0,\cdot)$ is the unique such polynomial (see~\cite{MR4637365}*{Section~2}) with
\begin{equation} \label{eq:Taylor_cond}
    \mbX^J f(x_0) = \mbX^J P_m[f](x_0,x) \, \big|_{x=x_0} \qquad \text{for all } d(J)\le m,
\end{equation}
that is,
\begin{equation} \label{eq:Taylor}
    P_m[f](x_0,x) = \sum_{d(J) \le m} a_J(x_0)\,\frac{(x_0^{-1}x)^J}{J!},
\end{equation}
where, by~\cite{MR4637365}*{Eq.~(2.10)} following~\cite{vodpup}*{Eq.~(2.9)},
\begin{equation} \label{eq:taylor_coefficients}
    a_J(x_0) = \sum_{\substack{|I| \le |J| \\ d(I) = d(J)}} \beta_{JI} \mbX^I f(x_0),
\end{equation}
for universal constants $\beta_{JI}$, depending only on $\G$ and the adapted basis, with $\beta_{JJ} = 1$ and $\beta_{JI} = 0$ whenever $d(I) = d(J)$, $|I| = |J|$, and $I \ne J$.

The two pointwise assertions of the next theorem follow from the stratified Taylor inequality of~\cite{FS82}*{Theorem~1.42 and Corollary~1.44}; for this class see also~\cite{MR4637365}*{Theorem~2.12}, and in step-two groups~\cite{arena2010taylor}*{Theorem~1}. The uniform H\"older remainder in the last assertion is the form used in~\cref{sec:whitney_jets}.

\begin{theorem}\label{thm:taylor_estimate}
    Let $\G$ be Carnot and $f \in \mathbf{C}_\G^m(\Omega)$, with $m\ge 1$. Then for each $x_0 \in \Omega$,
    \[
      \limsup_{x \to x_0}
      \frac{\bigl|f(x)-P_{m-1}[f](x_0,x)\bigr|}{d(x,x_0)^m}
      < \infty,
      \qquad
      \lim_{x \to x_0}
      \frac{\bigl|f(x)-P_m[f](x_0,x)\bigr|}{d(x,x_0)^m}
      =0.
    \]
    If, in addition, $f\in\mathbf C^{m,\alpha}_\G(\Omega)$ with $0<\alpha\le1$, then for every compact $K\Subset\Omega$ there are $C,r_0>0$ such that 
    \[
    |f(x)-P_m[f](x_0,x)|\le C d(x,x_0)^{m+\alpha}, \qquad  \forall \, x_0\in K, \  \forall \, x\in\Omega  \ \text{ s.t. } \ d(x,x_0)<r_0.
    \]
\end{theorem}

\begin{proof}
For $x_0\overline B_r(e)\subset\Omega$ set
\[
\eta_m(x_0,r):=\max_{d(J)=m}\sup_{\|z\|\le r} \bigl|\mbX^Jf(x_0z)-\mbX^Jf(x_0) \bigr|.
\]
The stratified Taylor inequality~\cite{FS82}*{Theorem~1.42} (see
also~\cite{FS82}*{Corollary~1.44}) gives constants $b_m,C_m\ge1$,
depending only on $m$, the group, and the homogeneous norm, such that every
$g\in C^\infty(\G)$ whose derivatives $\mbX^Jg(x_0)$, $d(J)\le m$, vanish
satisfies
\begin{equation}\label{eq:taylor_smooth}
    |g(x_0h)|\le C_m\|h\|^m\max_{d(J)=m}\sup_{\|z\|\le b_m\|h\|}
    |\mbX^Jg(x_0z)|.
\end{equation}
Indeed, \cite{FS82}*{Corollary~1.44} bounds the remainder by the differences
$\mbX^Jg(x_0z)-\mbX^Jg(x_0)$, which reduce to $\mbX^Jg(x_0z)$ when the jet of
$g$ at $x_0$ vanishes. Moreover, whatever family of left-invariant operators of
homogeneous degree $m$ is used for the supremum, on smooth functions each of
them is, by the Poincar\'e--Birkhoff--Witt theorem and homogeneity, a linear
combination of the ordered operators $\mbX^J$ with $d(J)=m$, which only changes the constant
$C_m$.

We pass to $\mathbf C^m_\G(\Omega)$ by mollification. Let
$x_0\overline B_{b_m\|h\|}(e)\subset\Omega$, fix $\chi\in C_c^\infty(\Omega)$
equal to $1$ near this compact set, and, for a mollifier
$\rho_\eps\in C_c^\infty(B_\eps(e))$ with
$\eps$ small, consider $f_\eps(x):=\int\rho_\eps(y)f(y^{-1}x)\,dy$. By unimodularity
of $\G$, $f_\eps(x)=\int\rho_\eps(xz^{-1})f(z)\,dz$, which is smooth on a
neighborhood of $\operatorname{supp}(\chi)$ for $\eps$ small. By left
invariance,
\[
\mbX^Jf_\eps(x)=\int\rho_\eps(y)\,\mbX^Jf(y^{-1}x)\,dy \quad \text{for
$d(J)\le m$} \implies \mbX^Jf_\eps\to\mbX^Jf
\]
locally uniformly. Apply~\eqref{eq:taylor_smooth} to
$g_\eps:=\chi f_\eps-P_m[\chi f_\eps](x_0,\cdot)$. By~\eqref{eq:Taylor_cond},
$\mbX^Jg_\eps(x_0)=0$ for $d(J)\le m$. For $d(J)=m$, homogeneity makes
$\mbX^JP_m[\chi f_\eps](x_0,\cdot)$ constant, equal to $\mbX^Jf_\eps(x_0)$, so
$\mbX^Jg_\eps(x_0z)=\mbX^Jf_\eps(x_0z)-\mbX^Jf_\eps(x_0)$
for $\|z\|\le b_m\|h\|$. Letting $\eps\to0$, with~\eqref{eq:taylor_coefficients}
for the coefficients, gives
\begin{equation}\label{eq:taylor_modulus_remainder}
    \bigl|f(x_0h)-P_m[f](x_0,x_0h)\bigr|
    \le C_m\|h\|^m\eta_m(x_0,b_m\|h\|)
    \quad\text{if }x_0\overline B_{b_m\|h\|}(e)\subset\Omega.
\end{equation}
Continuity of the top-order derivatives makes the right-hand side
of~\eqref{eq:taylor_modulus_remainder} $o(\|h\|^m)$ at each $x_0$, which is the second limit. By~\eqref{eq:taylor_coefficients},
$P_m[f](x_0,\cdot)-P_{m-1}[f](x_0,\cdot)$ is homogeneous of degree $m$ in
$h=x_0^{-1}x$, hence $O(\|h\|^m)$; with~\eqref{eq:taylor_modulus_remainder}
this gives the first assertion. Finally, if $K\Subset\Omega$ and
$K\overline B_{b_mr_0}(e)\Subset\Omega$, the finitely many $\mbX^Jf$,
$d(J)=m$, have a common $\alpha$-H\"older constant on this compact
thickening, so $\eta_m(x_0,b_m\|h\|)\le C\|h\|^\alpha$ for $x_0\in K$ and
$\|h\|<r_0$, and~\eqref{eq:taylor_modulus_remainder} gives the H\"older
estimate.
\end{proof}

Throughout this paper we use Pansu differentiability in the sense of~\cite{Pan89}. Let $\G_1,\G_2$ be graded groups with dilations $\delta^1_\lambda,\delta^2_\lambda$ and homogeneous distances $d_1,d_2$.

\begin{definition}\label{def:h-homomorphism}
A group homomorphism $L : \G_1 \to \G_2$ is an \emph{h-homomorphism} if $L\circ\delta^1_\lambda=\delta^2_\lambda\circ L$ for all $\lambda>0$. When $\G_1,\G_2$ are Carnot, a map $f : \Omega \to \G_2$, $\Omega\subset\G_1$ open, is \emph{Pansu differentiable} at $x \in \Omega$ if there is an h-homomorphism $d_H f(x) : \G_1 \to \G_2$ with
\[
d_2\bigl(f(x)^{-1} f(x'),\, d_H f(x)(x^{-1} x')\bigr)=o\bigl(d_1(x,x')\bigr) \qquad \text{as } x'\to x.
\]
\end{definition}

Every $F\in[\mathbf{C}_\G^1(\Omega)]^c$ is Pansu differentiable at every $p \in \Omega$, with differential $d_HF(p)(q)=\nabla_HF(p)\,(q_1,\ldots,q_\ell)^T$, where $q_i$ are the exponential coordinates of $q$. This follows componentwise from~\cref{thm:taylor_estimate} with $m=1$, since the displayed linear map on the first layer is an h-homomorphism by~\eqref{eq:BCH}.

%%%%%%%%%%%%%%%%%%%%%%%%%%%%%%%%%%%%%%%%%%%%%%%%%%%%%%%%%%%%%%%%%%%%
\section{Intrinsic rectifiability} \label{sec:rect_tech}
%%%%%%%%%%%%%%%%%%%%%%%%%%%%%%%%%%%%%%%%%%%%%%%%%%%%%%%%%%%%%%%%%%%%

\subsection{Intrinsic regular surfaces}

We start with the main definitions.

\begin{definition}\label{def:regular_surface}
    Let $1 \le k \le N$, $m \in \mathbb{N}$, and $\alpha \in [0,1]$. A set $S \subset \G$ is a \emph{$k$-dimensional $\G$-regular surface of class $(m,\alpha)$} (\emph{of class $m$} if $\alpha=0$) if every $p \in S$ has an open neighborhood $\Omega$ with a function $f \in [\mathbf{C}_\G^{m,\alpha}(\Omega)]^{N-k}$ such that:
    \begin{enumerate}[label=(\alph*)]
        \item for every $q\in\Omega$, the Pansu differential $d_Hf(q):\G\to\R^{N-k}$ is surjective and its kernel admits a homogeneous complement in $\G$;
        \item $S\cap\Omega=\{q\in\Omega \ : \ f(q)=0\}$.
    \end{enumerate}
\end{definition}

With the convention $\R^0=\{0\}$, the unique map to $\R^0$ is surjective with kernel $\G$, so the $N$-dimensional regular surfaces of any class are precisely the open subsets of $\G$. For class $1$, this notion was introduced in Heisenberg groups in~\cite{FSSC01} in codimension one and in~\cite{FSSC6} in low codimension, and in Carnot groups in~\cite{FSSC4} for hypersurfaces; see also~\cite{MSSC10}*{Definition~2.10} and, in low codimension in Carnot groups,~\cites{FranchiSerapioni2016IntrLip,AGV2022}.

\begin{definition}\label{def:tangent_group}
    Let $S$ be a $k$-dimensional $\G$-regular surface of class $1$, $p\in S$, and $f$ a local defining function. The \emph{homogeneous tangent subgroup} of $S$ at $p$ is 
    \[
    T_pS := \{ q\in\G \ : \ d_Hf(p)(q)=0\}.
    \]
\end{definition}

Notice that $T_pS$ does not depend on $f$~\cite{AGV2022}*{Lemma~2.14(iii)} and belongs to $\cG(\G,k)$ by~\cref{def:regular_surface}(a). Since $d_Hf(p)$ takes values in the first-layer group $\R^{N-k}$, it vanishes on every higher layer; thus $T_pS$ is vertical and normal, so $T_pS\in\cV(\G,k)$. By~\cref{prop:vertical_splitting}, every complement of $T_pS$ is horizontal; such a complement $\bV_p$ is generally non-canonical, and $\G=T_pS\rtimes\bV_p$. In the terminology of~\cite{AGV2022}, $p$ is a split-regular point of $f$.

\begin{definition}\label{def:rec}
    Let $1 \le k \le N$, $m\in\mathbb{N}$, and $\alpha\in[0,1]$. An $\cH^{k_{\m}}$-measurable set $E\subset\G$ is \emph{$(k,\G)$-rectifiable of class $(m,\alpha)$} if there are countably many $k$-dimensional $\G$-regular surfaces $S_i$ of class $(m,\alpha)$ with $\cH^{k_{\m}}(E\setminus\bigcup_{i}S_i)=0$.
\end{definition}

For $m=1$ and $\alpha=0$ this is rectifiability by regular surfaces as in the low-codimensional Heisenberg setting of~\cites{FSSC01,FSSC6,MSSC10}. In addition, \cite{FSSC03} proves codimension-one rectifiability of reduced boundaries in step-two Carnot groups.

\begin{definition}\label{def:intrinsic_graph}
    Let $1\le k\le N$, $\bW\in\cV(\G,k)$, and $\G=\bW\rtimes\bV$ with $\bV$ horizontal. A set $\Sigma\subset\G$ is an \emph{intrinsic graph} over $(\bW,\bV)$ if $\Sigma = \{ w\cdot\phi(w) \ : \ w\in U\}$ for some $U\subset\bW$ and $\phi:U\to\bV$. For $m\in\N$ and $\alpha\in[0,1]$, we write $\phi\in\mathbf{C}^{m,\alpha}_{\bW,\bV}(U)$ if $U$ is open in $\bW$, $\Sigma$ is a $k$-dimensional $\G$-regular surface of class $(m,\alpha)$, and
    \[
    T_{w\cdot\phi(w)}\Sigma\cap\bV=\{0\}\qquad\text{for every }w\in U.
    \]
    We call such a graph patch \emph{split-regular} and set $\mathbf C^m_{\bW,\bV}:=\mathbf C^{m,0}_{\bW,\bV}$.
\end{definition}

An \emph{intrinsic $\mathbf C_H^1$ graph} is one whose graph function is of class $\mathbf C^1_{\bW,\bV}$, in agreement with~\cite{AGV2022}*{Definition~2.13}. The notation $\mathbf C^{m,\alpha}_{\bW,\bV}$ constrains the graph, not the graph function since we never require $\phi$ to be differentiable along the left-invariant fields of $\bW$, which may fail. For $\phi:U\subset\bW\to\bV$ we write $\Phi_\phi(w):=w\cdot \phi(w)$, so that $\Phi_\phi(U)$ is the graph patch. Transversality to $\bV$ is open, so the intrinsic implicit-function theorem gives split-regular patches after shrinking. The following local representation is stated in~\cite{AGV2022}*{Lemma~2.14(i)}, whose proof follows from~\cite{AGV2022}*{Lemmas~2.10 and~2.12}.

\begin{theorem}\label{thm:intrinsic_surface_c1}
    Let $S\subset\G$ be a $k$-dimensional $\G$-regular surface of class $1$, $p\in S$, and $\G=\bW\rtimes\bV$ with $\bW\in\cV(\G,k)$ and $T_pS\cap\bV=\{0\}$. Then there are an open neighborhood $\Omega$ of $p$, an open set $U\subset\bW$ containing $\pi_{\bW}(p)$, and a continuous $\varphi\in\mathbf{C}^{1}_{\bW,\bV}(U)$ such that 
    \[
    S\cap\Omega = \bigl\{ w\cdot\varphi(w) \ : \ w\in U \bigr\}.
    \]
    In particular, $T_{w\cdot\varphi(w)}S\cap\bV=\{0\}$ for every $w\in U$.
\end{theorem}

The subgroup $\bV$ is horizontal by~\cref{prop:vertical_splitting}, applied
to $\bW$. In~\cite{AGV2022}
the hypothesis is that $T_pS$ is an intrinsic graph over $\bW$ along $\bV$;
this is equivalent to $T_pS\cap\bV=\{0\}$. Indeed, if
$\pi_\bW(t_1)=\pi_\bW(t_2)$ with $t_1,t_2\in T_pS$, then
$t_1^{-1}t_2=\pi_\bV(t_1)^{-1}\pi_\bV(t_2)\in T_pS\cap\bV$; so
$\pi_\bW|_{T_pS}$ is injective, hence open by invariance of domain, and its
image is a dilation-invariant open neighborhood of $e$ in $\bW$, hence
all of $\bW$.

%%%%%%%%%%%%%%%%%%%%%%%%%%%%%%%%%%%%%%%%%%%%%%%%%%%%%%%%%%%%%%%%%%%%
\subsection{Approximate tangent subgroups}\label{subsec:geometric_lemmas}
%%%%%%%%%%%%%%%%%%%%%%%%%%%%%%%%%%%%%%%%%%%%%%%%%%%%%%%%%%%%%%%%%%%%

\begin{definition}\label{def:tangent_space_star}
Following~\cite{San19}, let $s>0$, let $\mu$ be a Radon measure on $\G$, and let $x\in\G$. The \emph{$s$-dimensional approximate upper tangent cone} of $\mu$ at $x$ is
\[
    \Tan^{*s}(\mu,x) := \bigl\{v\in\G \ : \ \Theta^{*s}\bigl(\mu\res\mathcal X(x,v,\eps),x\bigr)>0
    \text{ for every }\eps>0\bigr\},
\]
where $\mathcal X(x,v,\eps) := \{y\in\G \ : \ d(v,\delta_r(x^{-1}y))<\eps \text{ for some } r>0\}$.
These directional cones are open, being unions of inverse images of open balls under continuous maps.
\end{definition}

For $x\in\G$, $\bW\in\cV(\G,k)$ with horizontal complement $\bV$, and $\eta>0$, the \emph{cone} with vertex $x$, axis $\bW$, and aperture $\eta$ is
\begin{equation}\label{eq:cone}
    \cC(x,\bW,\bV,\eta) := \bigl\{z\in\G \ : \ \|\pi_{\bV}(x^{-1}z)\| \le \eta\,\|\pi_{\bW}(x^{-1}z)\| \bigr\}.
\end{equation}
The following result is the homogeneous-group version of~\cite{San19}*{Lemma~3.6}. We give the proof for completeness.

\begin{lemma}\label{lem:tangent_condition}
    Fix $1\le k\le N$. Let $A\subset\G$ be $\cH^{k_{\m}}$-measurable and locally $\cH^{k_{\m}}$-finite, $x\in\G$, and $\G=\bW\rtimes\bV$ with $\bW\in\cV(\G,k)$. The following are equivalent:
    \begin{enumerate}[label=\normalfont(\roman*)]
        \item $\Tan^{*k_{\m}}(\cH^{k_{\m}}\res A,x)\subseteq\bW$;
        \item $\Theta^{*k_{\m}}\bigl(\cH^{k_{\m}}\res(A\setminus\cC(x,\bW,\bV,\eta)),x\bigr)=0$ for every $\eta>0$;
        \item for every $\eps>0$,
        \[
        \lim_{r\downarrow0}r^{-k_{\m}}
        \cH^{k_{\m}}\Bigl(
        A\cap B_r(x)
        \cap
        \bigl\{z\in\G \ : \ 
        \|\pi_{\bV}(x^{-1}z)\|>\eps r
        \bigr\}
        \Bigr)=0.
        \]
    \end{enumerate}
\end{lemma}

\begin{proof}
Since $k\ge1$, $s:=k_{\m}>0$, and $\bV$ is horizontal by~\cref{prop:vertical_splitting}. Write $\mu:=\cH^s\res A$, a locally finite Radon measure, $h(g):=\|\pi_{\bW}(g)\|$, and
$v(g):=\|\pi_{\bV}(g)\|$ (both functions are continuous and one-homogeneous). For
$\eta>0$ let $D_\eta:=\{g \ : \ v(g)>\eta h(g)\}$, so that
$A\setminus\cC(x,\bW,\bV,\eta)=A\cap xD_\eta$.

\smallskip
\textup{(i)$\iff$(ii).} If $K\subset\G\setminus\Tan^{*s}(\mu,x)$ is
compact, then
\begin{equation}\label{eq:compact_directions_zero_density}
    \Theta^{*s}\bigl(\mu\res\mathcal X(x,K),x\bigr)=0,
    \quad
    \mathcal X(x,K):=\{y \ : \ \delta_r(x^{-1}y)\in K
    \text{ for some }r>0\}.
\end{equation}
Here $\mathcal X(x,K)=\bigcup_{n\in\N}x\{\delta_tu \ : \ u\in K,\ 1/n\le t\le n\}$ is a countable union of compact sets, hence Borel.
Each $u\in K$ has $\eps_u>0$ with
$\Theta^{*s}(\mu\res\mathcal X(x,u,\eps_u),x)=0$, and finitely many balls
$B_{\eps_u}(u)$ cover $K$, so $\mathcal X(x,K)$ lies in a finite union of
such cones.
Assume~\textup{(i)} and fix $\eta>0$. The compact set
$K_\eta:=\{\|g\|=1,\ v(g)\ge\eta h(g)\}$ is disjoint from $\bW$, hence from
$\Tan^{*s}(\mu,x)$, and $D_\eta\subset\{\delta_tu \ : \ t>0,\ u\in K_\eta\}$,
so~\eqref{eq:compact_directions_zero_density} gives~\textup{(ii)}.
Conversely, assume~\textup{(ii)} and let $u\notin\bW$. Then $v(u)>0$, and by
continuity and homogeneity there are $\eta,\eps>0$ with
$\mathcal X(x,u,\eps)\subset xD_\eta$; thus
$\Theta^{*s}(\mu\res\mathcal X(x,u,\eps),x)=0$ and
$u\notin\Tan^{*s}(\mu,x)$.

\textup{(ii)$\implies$(iii).} If $z\in B_r(x)$ and $v(x^{-1}z)>\eps r$,
then $v(x^{-1}z)>(\eps/C_{\bW,\bV})h(x^{-1}z)$ by~\eqref{eq:piVcontin}, so
the set in~\textup{(iii)} lies in $A\cap xD_{\eps/C_{\bW,\bV}}$.

\textup{(iii)$\implies$(ii).} Fix $\eta>0$, choose $\eps>0$ with
$\tau:=\eps(1+\eta^{-1})<1$, and set $F(r):=\mu(B_r(x)\cap xD_\eta)$ and
$E_\eps(r):=A\cap B_r(x)\cap\{z \ : \ v(x^{-1}z)>\eps r\}$. If
$z\in A\cap B_r(x)\cap xD_\eta$ but $z\notin E_\eps(r)$, then
$h(x^{-1}z)<\eta^{-1}v(x^{-1}z)\le\eta^{-1}\eps r$, hence
$\|x^{-1}z\|\le h(x^{-1}z)+v(x^{-1}z)<\tau r$. Thus
$F(r)\le\mu(E_\eps(r))+F(\tau r)$, and iterating,
\[
    \frac{F(r)}{r^s}
    \le
    \sum_{j=0}^{\infty}\tau^{js}
    \frac{\mu(E_\eps(\tau^jr))}{(\tau^jr)^s}
    \le
    \frac{1}{1-\tau^s}
    \sup_{0<t\le r}\frac{\mu(E_\eps(t))}{t^s}.
\]
For fixed $r$, the sets $B_{\tau^nr}(x)\cap xD_\eta$ decrease to the empty set because $x\notin xD_\eta$. Local finiteness therefore gives $F(\tau^nr)/r^s\to0$, justifying the iteration, and $s>0$ ensures convergence of the geometric series. The right-hand side tends to zero by~\textup{(iii)}.
\end{proof}

Condition~\textup{(i)} does not involve $\bV$. Hence, if \textup{(ii)}
or~\textup{(iii)} holds for one horizontal complement of $\bW$, it holds for
every other one; this is used tacitly when complements are changed in the
proofs of~\cref{thm:c1_rec,thm:main}.

%%%%%%%%%%%%%%%%%%%%%%%%%%%%%%%%%%%%%%%%%%%%%%%%%%%%%%%%%%%%%%%%%%%%%%%%%%%%%%%%%%%%%%%%%%%%%%%%%%%%%%%%%%%%%%%%%%%%%%%%%%%%
\section{Whitney extension and adapted graph jets} \label{sec:whitney_jets}
%%%%%%%%%%%%%%%%%%%%%%%%%%%%%%%%%%%%%%%%%%%%%%%%%%%%%%%%%%%%%%%%%%%%%%%%%%%%%%%%%%%%%%%%%%%%%%%%%%%%%%%%%%%%%%%%%%%%%%%%%%%%

This section collects the tools for the proofs of~\cref{thm:c1_rec,thm:main}.
More precisely, \cref{subsec:whitney} recalls the Whitney extension theorems, which turn
compatible jets on a closed set into ambient $\mathbf C^{m,\alpha}_\G$
functions, whose non-critical level sets are regular surfaces
(\cref{lem:level_set}). \Cref{subsec:graph_jets} works on a fixed first-order graph patch
$M$. To be precise, \cref{lem:surface_normal_form_jet} converts a regular surface
through a point of $M$ into a normal-form polynomial for the split
projection; \cref{lem:graph_blowup} identifies the blow-ups of $M$;
\cref{lem:normal_form_norming} bounds normal-form
polynomials by their values on subsets of $M$ of positive density; and
\cref{lem:trace_normalized_selection} selects
compatible jets and builds regular surfaces from them. These results are used
in Steps~2--4 of the proof of~\cref{thm:higher_order}.

\subsection{Whitney extension}\label{subsec:whitney}

For a real-valued jet $\{f_J\}_{d(J)\le m}$ on a closed set $F\subset\G$ and $x_0\in F$, let $P_m[\{f_J\}_{d(J)\le m}](x_0,\cdot)$ be the polynomial~\eqref{eq:Taylor} whose coefficients are given by~\eqref{eq:taylor_coefficients} with each $\mbX^If(x_0)$ replaced by $f_I(x_0)$. The linear map sending a polynomial of homogeneous degree at most $m$ centered at $x_0$ to its jet $(\mbX^JP(x_0))_{d(J)\le m}$ is injective by the uniqueness in~\eqref{eq:Taylor_cond}, hence bijective, both spaces having dimension $\#\{J:d(J)\le m\}$. By the same uniqueness every such polynomial is its own Taylor polynomial (see also~\cite{MR4637365}*{Proposition~2.13}), so~\eqref{eq:taylor_coefficients} is the inverse map. Therefore
\begin{equation}\label{eq:jet_polynomial_interpolation}
    \mbX^J P_m\bigl[\{f_I\}_{d(I)\le m}\bigr](x_0,x)\big|_{x=x_0}
    =f_J(x_0),
    \qquad d(J)\le m.
\end{equation}
We use the Whitney extension theorem of Vodopyanov--Pupyshev~\cite{vodpup}*{Theorem~4}, in the formulation of~\cite{MR4637365}*{Theorem~2.17}.

\begin{theorem}\label{thm:whitney_extension}
Let $m\in\N$, $F\subset\G$ closed, and $\{f_J\}_{d(J)\le m}$ a real-valued jet on $F$ such that \textup{(i)} $\sup_{K}|f_J|<\infty$ for every compact $K\subset F$ and every $d(J)\le m$, and \textup{(ii)} for every $\bar x\in F$ and $\eps>0$ there is $\delta>0$ such that
\begin{equation}\label{eq:whitney_little_o_condition}
\left|f_J(y)-\mbX^J P_m\bigl[\{f_I\}_{d(I)\le m}\bigr](x,y)\right|
\le \eps\, d(x,y)^{m-d(J)}
\end{equation}
whenever $x\ne y\in F\cap B_\delta(\bar x)$ and $d(J)\le m$. Then there is $u\in\mathbf C^m_\G(\G)$ with $\mbX^J u|_F=f_J$ for every $d(J)\le m$.
\end{theorem}

\begin{definition}\label{def:lip_jet}
Let $m\in\N$ and $m<\gamma\le m+1$. A jet $\{f_J\}_{d(J)\le m}$ on a closed set $F\subset\G$ belongs to $\operatorname{Lip}(\gamma,F)$ if there is $C_F>0$ such that, for all $x,y\in F$ and $d(J)\le m$,
\begin{equation}\label{eq:lip_jet_condition}
    |f_J(x)|\le C_F,
    \qquad
    \left| f_J(y)-\mbX^J P_m\bigl[\{f_I\}_{d(I)\le m}\bigr](x,y)\right|
        \le C_F\, d(x,y)^{\gamma-d(J)}.
\end{equation}
We write $u\in\operatorname{Lip}(\gamma,\G)$ if the derivatives $\mbX^Ju$, $d(J)\le m$, exist and their jet satisfies~\eqref{eq:lip_jet_condition} on the whole group $F=\G$.
\end{definition}

\begin{theorem}\label{thm:whitney_alpha_holder}
Let $m\in\N$, $\alpha\in(0,1]$, and $F\subset\G$ closed. If $\{f_J\}_{d(J)\le m}\in\operatorname{Lip}(m+\alpha,F)$, there is $u\in\operatorname{Lip}(m+\alpha,\G)$ with $\mbX^J u|_F=f_J$ for all $d(J)\le m$. Precisely, $u\in\mathbf C^{m,\alpha}_\G(\G)$.
\end{theorem}

This is~\cite{vodpup}*{Theorem~2 and Lemma~6}, as stated in~\cite{MR4637365}*{Theorem~2.21}. The final inclusion is~\cite{MR4637365}*{Lemma~2.20}. Note that by~\eqref{eq:jet_polynomial_interpolation} the polynomial term in~\eqref{eq:lip_jet_condition} tends to $\mbX^Ju(x)$ as $y\to x$, so every $\mbX^Ju$ is continuous. If $d(J)=m$, this term is homogeneous of degree zero in $x^{-1}y$, hence equal to $\mbX^Ju(x)$, and~\eqref{eq:lip_jet_condition} is the $\alpha$-H\"older estimate.

\begin{lemma}\label{lem:level_set}
Let $m\in\mathbb N$ and $\alpha\in[0,1]$, let $\bV$ be a horizontal subgroup of
dimension $1\le c<N$, let $V_1,\ldots,V_c$ be a basis of
$\Lie(\bV)$, and let $F\in[\mathbf C^{m,\alpha}_\G(\G)]^c$. If $F=0$ and
$V_iF=e_i$, for $i=1,\ldots,c$, on a set $K\subset\G$, then $K$ is contained in an
$(N-c)$-dimensional $\G$-regular surface of class $(m,\alpha)$.
\end{lemma}

\begin{proof}
The matrix $(V_jF_i(q))_{i,j}$ represents $d_HF(q)|_\bV$ in the given basis,
so 
\[
\Omega:= \left\{q \ : \ \|(V_jF_i(q))_{i,j}-I_c\|_{\mathrm{op}}<\tfrac12\right\}
\]
is open, contains $K$, and $d_HF(q)|_\bV$ is invertible for $q\in\Omega$. For such
$q$, $d_HF(q)$ is surjective and $\ker d_HF(q)\cap\bV=\{0\}$. Moreover, for
$g\in\G$ there is $v\in\bV$ with $d_HF(q)(v)=d_HF(q)(g)$, so
$gv^{-1}\in\ker d_HF(q)$. Hence the normal homogeneous subgroup
$\ker d_HF(q)$ is complemented by $\bV$, and by~\cref{def:regular_surface}
the set $F^{-1}(0)\cap\Omega$ is an $(N-c)$-dimensional $\G$-regular surface
of class $(m,\alpha)$ containing $K$.
\end{proof}

%%%%%%%%%%%%%%%%%%%%%%%%%%%%%%%%%%%%%%%%%%%%%%%%%%%%%%%%%%%%%%%%%%%%%%%%%%%%%%%%
\subsection{Adapted graph jets and regular selection}\label{subsec:graph_jets}
%%%%%%%%%%%%%%%%%%%%%%%%%%%%%%%%%%%%%%%%%%%%%%%%%%%%%%%%%%%%%%%%%%%%%%%%%%%%%%%%

Since the parameter subgroups $\bW$ need not be Carnot
(\cref{rmk:graded_not_carnot}), Carnot-group Whitney and Lusin theorems do
not apply on them, and we thus work with ambient jets on $\G$. Throughout this
subsection we fix $1\le k\le N$, a subgroup
$\bW\in\cV(\G,k)$, and a horizontal subgroup $\bV$ with the splitting $\G=\bW\rtimes\bV$.
We set $c:=\dim\bV=N-k$ and $s:=k_{\m}$, and identify $\bV$ with $\R^c$ through a
basis $V_1,\ldots,V_c$ of $\Lie(\bV)$. We write $\Pi_\bV:\G\to\R^c$ for the
coordinates of the split projection onto $\bV$, which is a homogeneous
homomorphism equal to the identity on $\bV$, so $V_i\Pi_\bV=e_i$. Finally,
$M=\Phi_\phi(U)$ is a split-regular $\mathbf C^1_H$ intrinsic graph over an open
set $U\subset\bW$, with values in $\bV$. We take $M$ relatively compact,
with its closure contained in a larger split-regular graph chart. Thus
$\bV$ complements $T_xM$ for every $x\in M$, and $\cH^s(M)<\infty$
by~\cite{AGV2022}*{Proposition~2.17}.

Constants attached to a point, such as $C_p,r_p,\omega_p$
in~\cref{lem:surface_normal_form_jet} and $L_p$
in~\cref{lem:trace_normalized_selection}, are only pointwise, while uniformity on
compact pieces is gained only in Step~3 of the proof
of~\cref{lem:trace_normalized_selection}, through~\cref{lem:compact_pieces}.

We denote by $\mathcal N_{\bV}^m$ the finite-dimensional space of scalar
polynomials $P$ of homogeneous degree at most $m$ with $P(xv)=P(x)$ for all
$x\in\G$ and $v\in\bV$ or, equivalently, such that $P=P\circ\pi_{\mathbb T}$ for any
splitting $\G=\mathbb T\rtimes\bV$. This space depends neither on a base
point nor on a complement of $\bV$, and $V_iP\equiv0$ for every
$P\in\mathcal N_\bV^m$.

\begin{lemma}\label{lem:surface_normal_form_jet}
Let $m\in\mathbb N$ and $\alpha\in[0,1]$, let $S$ be a $k$-dimensional $\G$-regular surface of class $(m,\alpha)$, and let $p\in S$ with $T_pS\cap\bV=\{0\}$, so that $\G=T_pS\rtimes\bV$. Then there are $r_p>0$, $C_p<\infty$ if $\alpha>0$, and $R_p\in[\mathcal N_\bV^m]^c$ with $R_p(p)=\Pi_\bV(p)$ such that
\begin{equation}\label{eq:surface_normal_form_estimate}
    |\Pi_{\bV}(q)-R_p(q)| \le
    \begin{cases}
        C_p\, d(p,q)^{m+\alpha}, & \alpha>0,\\
        \omega_p(d(p,q))\,d(p,q)^{m}, & \alpha=0,
    \end{cases}
\end{equation}
holds for every $q\in S\cap B_{r_p}(p)$, where $\omega_p(r)\downarrow0$ as $r\downarrow0$.
\end{lemma}

\begin{proof}
If $c=0$, there is nothing to prove. Let $c>0$, and let
$F\in[\mathbf C^{m,\alpha}_{\G}(\Omega)]^{c}$ define $S$ in an open
neighborhood $\Omega$ of $p$. Since
$\ker d_HF(p)=T_pS$ is complemented by $\bV$, the map $d_HF(p)|_{\bV}$ is
invertible so, replacing $F$ by $(d_HF(p)|_\bV)^{-1}\circ F$, we may assume
$d_HF(p)|_{\bV}=\operatorname{Id}$. Let $\varpi:\G\to\bV$ be the projection
of the splitting $T_pS\rtimes\bV$, and write $g=\tau(g)\cdot\varpi(g)$ with
$\tau(g)\in T_pS$. By~\eqref{eq:projection_algebra}, $\varpi$ is a
h-homomorphism, hence linear in exponential coordinates, and
by~\eqref{eq:BCH} $\tau$ is polynomial and commutes with dilations.

Set $P:=P_m[F](p,\cdot)$ and $\widehat P(\tau,v):=P(p\cdot\tau\cdot v)$ for
$\tau\in T_pS$ and $v\in\bV$, an $\R^c$-valued polynomial. We use the
product dilations $(\tau,v)\mapsto(\delta_t\tau,\delta_tv)$, for which $v$
has degree one. For them,
$\widehat P=\widehat P_1+\cdots+\widehat P_m$ splits into homogeneous
components with zero constant term because $F(p)=0$, and
$\widehat P_1(\tau,v)=v$ by the normalization of $d_HF(p)$. Writing $[Q]_r$
for the component of degree $r$ of a polynomial $Q$, set $\psi_1:=0$ and, for
$2\le r\le m$, define the $\R^c\cong\bV$-valued polynomial $\psi_r$ on
$T_pS$, homogeneous of degree $r$, by
\begin{equation}\label{eq:normal_form_homogeneous_recursion}
    \psi_r(\tau):=-\bigl[\widehat P(\tau,\psi_{<r}(\tau))\bigr]_r,
    \qquad \psi_{<r}:=\sum_{h<r}\psi_h.
\end{equation}
Every monomial of $\widehat P-\widehat P_1$ has homogeneous degree at least two.
Replacing a degree-one $v$-factor by $\psi_r$ therefore produces degree at
least $r+1$. Thus replacing $\psi_{<r}$ by $\psi_{<r}+\psi_r$ changes degrees
at most $r$ only through $\widehat P_1=v$, and the degree-$r$ component becomes
$\psi_r+[\widehat P(\tau,\psi_{<r})]_r=0$. Hence
$\psi_p:=\sum_{r=2}^m\psi_r$ ($\psi_p:=0$ if $m=1$) satisfies $\psi_p(e)=0$
and
\begin{equation}\label{eq:normal_form_formal_remainder}
    P\bigl(p\cdot\tau\cdot\psi_p(\tau)\bigr)
    =O(\|\tau\|^{m+1}) \qquad\text{as }\tau\to e.
\end{equation}
We can now define the candidate polynomial as
\begin{equation}\label{eq:normal_form_projection_identity}
    R_p(x):=\Pi_{\bV}(p)+\Pi_{\bV}(\tau(g))+\psi_p(\tau(g)),
    \qquad g:=p^{-1}x.
\end{equation}
For $m=1$, $R_p(x)=\Pi_\bV(x)-\varpi(p^{-1}x)$. In general, $R_p$ is a polynomial of homogeneous degree at most $m$ in $g$, hence in $x$
by~\eqref{eq:BCH}. Since $\tau(gv')=\tau(g)$ for $v'\in\bV$, it is constant on
left cosets of $\bV$, so $R_p\in[\mathcal N_\bV^m]^c$. Moreover,
$R_p(p)=\Pi_\bV(p)$ because $\psi_p(e)=0$.

Let $q\in S\cap B_{r_p}(p)$, where $r_p>0$ is fixed below with
$B_{r_p}(p)\subset\Omega$ and $r_p$ at most the radius $r_0$
of~\cref{thm:taylor_estimate} for $K=\{p\}$, and set $\tau:=\tau(p^{-1}q)$ and
$v:=\varpi(p^{-1}q)$. Since $\Pi_\bV$ is a homomorphism equal to the identity
on $\bV$, we have 
\[
\Pi_{\bV}(q)=\Pi_{\bV}(p)+\Pi_{\bV}(\tau)+v \implies \Pi_\bV(q)-R_p(q)=v-\psi_p(\tau).
\]
Since $F(q)=0$,~\cref{thm:taylor_estimate}
gives $|P(q)|\le Cd(p,q)^{m+\alpha}$ if $\alpha>0$ and $|P(q)|=o(d(p,q)^m)$ if
$\alpha=0$, while~\eqref{eq:piVcontin} and~\eqref{eq:normal_form_formal_remainder}
give 
\[
|P(p\cdot\tau\cdot\psi_p(\tau))|\le Cd(p,q)^{m+1}.
\]
The maps
$\Gamma_\tau(w):=P(p\cdot\tau\cdot w)$, $w\in\bV$, satisfy
$D_w\Gamma_e(0)=\operatorname{Id}$, so $\|D_w\Gamma_\tau-\operatorname{Id}\|_{\mathrm{op}}\le\frac12$
for $\tau$ near $e$ and $w$ in a Euclidean ball about $0$. We choose $r_p$ so
small that this ball contains $v$ and $\psi_p(\tau)$ for every
$q\in S\cap B_{r_p}(p)$. Integrating along segments gives
\[
    |v-\psi_p(\tau)|
    \le2\bigl|\Gamma_\tau(v)-\Gamma_\tau(\psi_p(\tau))\bigr|
    \le2|P(q)|+2Cd(p,q)^{m+1},
\]
which proves~\eqref{eq:surface_normal_form_estimate}. The conclusion for $\alpha=0$ follows if one takes, for instance, $\omega_p(r):=\sup\bigl(\{0\}\cup\{|\Pi_\bV(q)-R_p(q)|\,d(p,q)^{-m} \ : \ q\in S,\
0<d(p,q)\le r\}\bigr)$ for $r<r_p$.
\end{proof}

\begin{workedexample}
In $\bH^1$, use exponential coordinates of weights $(1,1,2)$, with product
\[
 (x,y,t)\cdot(x',y',t')
 =(x+x',y+y',t+t'+\tfrac12(xy'-yx')).
\]
Let $\bW=\{x=0\}$, $\bV=\{y=t=0\}$, $\Pi_{\bV}(x,y,t)=x$, and
$X=\partial_x-\frac y2\partial_t$. The surface $S=\{F=0\}$ defined by
$F(x,y,t)=x-y^2-t$, is regular near $0$ since $XF=1+y/2$, and has tangent space $T_0S=\bW$.
For $w=(0,y,s)\in\bW$,
\[
 F\bigl(w\cdot(v,0,0)\bigr)=(1+\tfrac y2)v-y^2-s \stackrel{\eqref{eq:normal_form_homogeneous_recursion}}{\implies}\psi_2(w)=y^2+s.
\]
Since $(x,y,t)=(0,y,t+xy/2)\cdot(x,0,0)$, the ambient polynomial is
$R_0(x,y,t)=y^2+t+\tfrac12xy$, of homogeneous degree two, constant on left
cosets of $\bV$, with $R_0(0)=0$ and $XR_0\equiv0$. For $q=(x,y,t)\in S$ near
$0$ we have $x=y^2+t=O(\|q\|^2)$ and $y=O(\|q\|)$, which implies the $(2,1)$ normal-form estimate:
\[
|\Pi_{\bV}(q)-R_0(q)|=\tfrac12|xy|=O(\|q\|^3).
\]
In addition, $X(\Pi_{\bV}-R_0)=1$ is the transversality required
in~\cref{lem:level_set}.
\end{workedexample}

\begin{lemma}\label{lem:graph_blowup}
For every $a\in M$,
\[
    \cH^s\res\delta_{1/r}(a^{-1}M)\overset{*}{\rightharpoonup}\cH^s\res T_aM
    \qquad\text{locally as }r\downarrow0.
\]
\end{lemma}

\begin{proof}
Let $F\in[\mathbf C^1_\G(\Omega)]^c$ define $M$ near $a$, so that
$M\cap\Omega'=\{F=0\}\cap\Omega'$ for a neighborhood $\Omega'\subset\Omega$ of $a$.
Since $\bV$ complements $\ker d_HF(a)=T_aM$, the point $a$ is split-regular
for $F$ in the sense of~\cite{AGV2022}. For the level offset $b=0$,
\cite{AGV2022}*{Lemma~3.4} states that the $s$-dimensional Hausdorff (or
spherical) measures restricted to $\delta_{1/r}(a^{-1}\{F=0\})$ converge weak*
to the same measure restricted to $\ker d_HF(a)=T_aM$. The normalization of
$\cH^s$ changes these measures by a fixed factor, and the dilates
$\delta_{1/r}(a^{-1}\Omega')$ exhaust $\G$, so compactly supported test functions do not distinguish $M$ from $\{F=0\}$ for $r \ll 1$.
\end{proof}

The next estimate converts bounds on a subset of the graph of positive measure
into bounds for all ambient derivatives of a coset-invariant polynomial.

\begin{lemma}\label{lem:normal_form_norming}
Let $K\subset M$ be compact and $m\in\mathbb N$. For every $\theta>0$ there
are $c_0,r_0>0$ with the following property. If $x\in K$, $0<r<r_0$,
$P\in[\mathcal N_{\bV}^m]^c$, and $Z\subset M\cap B_r(x)$ is Borel with
$\cH^s(Z)\ge\theta r^s$, then, with the Euclidean norm on $\R^c$,
\begin{equation}\label{eq:normal_form_norming_one_point}
    \sum_{d(J)\le m}r^{d(J)}|\mbX^JP(x)| \le c_0^{-1}\sup_{Z}|P|.
\end{equation}
\end{lemma}

\begin{proof}
If $c=0$, there is nothing to prove. Let $c>0$ and notice that it suffices to treat scalar $P$ since applying that case to each component and
summing gives the general case after replacing $c_0$ by $c_0/c$.
Let $\pi_\bW$ be the projection of the splitting $\G=\bW\rtimes\bV$, and write
$|\cdot|_\bW:=\cH^s\res\bW$, a Haar measure on $\bW$
by~\cref{rmk:haar_measure}. Since $g\in\pi_\bW(g)\bV$, every
$P\in\mathcal N_\bV^m$ satisfies $P(g)=P(\pi_\bW(g))$.

Since $K$ is compact and $\pi_\bW$ is continuous, there are an open set
$U'\Subset U$ and $r_0>0$ such that $\pi_\bW(M\cap B_{r_0}(x))\subset U'$ for
every $x\in K$. The intrinsic area formula~\cite{AGV2022}*{Theorem~1.1}, whose
area factor is continuous on $U$ by~\cite{AGV2022}*{Lemma~3.2}, gives
$J_0<\infty$ with $\cH^s(\Phi_\phi(B))\le J_0|B|_\bW$ for every Borel
$B\subset U'$.

The graph map $\Phi_\phi:U\to M$ is a homeomorphism with inverse
$\pi_\bW|_M$, so $\pi_\bW(Z)$ is Borel whenever $Z\subset M$ is Borel.
Fix $x$, $r$, $P$, and $Z$ as in the statement. Then $Z=\Phi_\phi(\pi_\bW(Z))$,
so $|\pi_\bW(Z)|_\bW\ge\theta r^s/J_0$. Write $x=w_xv$, with
$w_x:=\pi_\bW(x)$ and $v:=\pi_\bV(x)$. By~\eqref{eq:projection_algebra},
\[
    \pi_\bW(x^{-1}y)=v^{-1}\bigl(w_x^{-1}\pi_\bW(y)\bigr)v.
\]
Thus $\pi_\bW(x^{-1}Z)$ is obtained from $\pi_\bW(Z)$ by a left translation
in $\bW$ followed by conjugation by $v^{-1}$. In exponential coordinates this
conjugation is the linear map
$\operatorname{Ad}_{v^{-1}}|_{\mathfrak w}=\exp(-\operatorname{ad}_{\log v}|_{\mathfrak w})$,
$\mathfrak w:=\Lie(\bW)$, whose determinant is $1$ because
$\operatorname{ad}_{\log v}$ is nilpotent on the ideal $\mathfrak w$. Hence
$|\pi_\bW(x^{-1}Z)|_\bW=|\pi_\bW(Z)|_\bW$.

Set $\widetilde P(z):=P(x\cdot\delta_rz)$ and $Z':=\delta_{1/r}\pi_\bW(x^{-1}Z)$. Left
translations and dilations preserve the degree bound and map left cosets of $\bV$ to left cosets of $\bV$,
so $\widetilde P\in\mathcal N_\bV^m$. Homogeneity of the left-invariant fields gives
$\mbX^J\widetilde P(e)=r^{d(J)}\mbX^JP(x)$, and, by $\bV$-invariance of $P$,
$\sup_{Z'}|\widetilde P|=\sup_Z|P|$. Moreover, by~\eqref{eq:piVcontin}, $Z'$ lies in the
compact set $B':=\bW\cap\overline B_{C'}(e)$ with $C':=C_{\bW,\bV}$, and
$|Z'|_\bW=r^{-s}|\pi_\bW(x^{-1}Z)|_\bW\ge\kappa:=\theta/J_0$.

It remains to find $c_0>0$, depending only on $m$, $C'$, and $\kappa$, with
\[
    \sum_{d(J)\le m}|\mbX^J\widetilde P(e)|\le c_0^{-1}\sup_{Z'}|\widetilde P|
\]
for every $\widetilde P\in\mathcal N_\bV^m$ and every Borel set $Z'\subset B'$ with
$|Z'|_\bW\ge\kappa$.
On the finite-dimensional space $\mathcal N_\bV^m$, the left-hand side is a
norm by the uniqueness in~\eqref{eq:Taylor_cond}, and so is
$\widetilde P\mapsto\sup_{B'}|\widetilde P|$. Indeed, a polynomial vanishing on $B'$ vanishes on $\bW$,
because $B'$ has nonempty interior in $\bW$, and hence everywhere by
$\bV$-invariance. It therefore suffices to bound $\sup_{B'}|\widetilde P|$ by
$\sup_{Z'}|\widetilde P|$. We argue by contradiction:

If no such bound held, there would be $P_j\in\mathcal N_\bV^m$
and Borel sets $Z'_j\subset B'$ with $\sup_{B'}|P_j|=1$,
$|Z'_j|_\bW\ge\kappa$, and $\sup_{Z'_j}|P_j|\to0$. After passing to a
subsequence, $P_j\to P_\infty$ uniformly on $B'$, with $\sup_{B'}|P_\infty|=1$, and
$Z'_j\subset\{w\in B' \ : \ |P_\infty(w)|\le\eta_j\}$ for some $\eta_j\to0$. In exponential
coordinates, $P_\infty|_\bW$ is a nonzero polynomial, so its zero set is
$|\cdot|_\bW$-null and $|\{w\in B' \ : \ |P_\infty(w)|\le\eta\}|_\bW\to0$ as
$\eta\downarrow0$. This contradicts $|Z'_j|_\bW\ge\kappa$.
\end{proof}

Thus~\cref{lem:normal_form_norming} is a Remez-type inequality on $\bW$. The
graph structure of $M$ enters only through the area bound
$\cH^s(\Phi_\phi(B))\le J_0|B|_\bW$.

\begin{lemma}\label{lem:trace_normalized_selection}
Let $m\in\mathbb N$ and $\alpha\in[0,1]$, and let $A\subset M$ be
$\cH^s$-measurable with
\begin{equation}\label{eq:trace_selection_density_tangent}
    \Theta_*^{s}(A,p)>0
    \qquad\text{for $\cH^s$-a.e.\ $p\in A$.}
\end{equation}
Suppose that for $\cH^s$-a.e.\ $p\in A$ there is an $\R^c$-valued polynomial
$R_p$ (recall that $\bV\cong\R^c$) with the following properties:
\begin{enumerate}[label=\textup{(\roman*)}]
\item $R_p\in[\mathcal N_\bV^m]^c$;
\item $R_p(p)=\Pi_{\bV}(p)$;
\item if $\alpha>0$, then for some $L_p<\infty$,
\[
    \lim_{r\downarrow0}r^{-s}\cH^{s}\Bigl( A\cap B_r(p)\cap \bigl\{q \ : \ |\Pi_{\bV}(q)-R_p(q)|>L_p r^{m+\alpha}\bigr\} \Bigr)=0;
\]
if $\alpha=0$, then the same limit holds with $r^m$ and with every constant $L>0$.
\end{enumerate}
Then $A$ is covered, up to an $\cH^s$-null set, by countably many compact
sets $K\subset A$ on each of which the jet
\begin{equation}\label{eq:trace_selection_whitney_coeffs}
    A_J(p):=\mbX^J R_p(p), \qquad p\in K, \quad d(J)\le m,
\end{equation}
belongs componentwise to $\operatorname{Lip}(m+\alpha,K)$ if $\alpha>0$, and
satisfies componentwise the hypotheses of~\cref{thm:whitney_extension} if
$\alpha=0$, for a Borel choice of polynomials $R_p$ satisfying
\textup{(i)--(iii)} (unique by~\cref{rmk:jet_uniqueness}). Consequently, $A$ is covered, up to an $\cH^s$-null set, by countably many
$k$-dimensional $\G$-regular surfaces of class $(m,\alpha)$.
\end{lemma}

\begin{proof}
\proofstep{1 (Reductions)} Removing an $\cH^s$-null set changes neither
densities nor exceptional-set estimates. We may therefore assume that $A$ is
Borel and that every $p\in A$ satisfies~\eqref{eq:trace_selection_density_tangent}
and admits a polynomial satisfying \textup{(i)--(iii)}. Since $A\subset M$ and
$\cH^s(M)<\infty$, the measure $\mu:=\cH^s\res A$ is finite. If $c=0$, then $k=N$, the jets are
empty, inner regularity gives the compact sets, and $\G$ itself is a regular
surface of every class covering $A$. Consequently, we may also assume $c>0$.

\proofstep{2 (Borel-choice)} Fix a basis $B_1,\ldots,B_{D_m}$ of
$\mathcal N_\bV^m$, $D_m:=\dim\mathcal N_\bV^m$, and, for $\xi=(\xi_{bi})\in Y:=\R^{D_m\times c}$, let
$R^\xi:=(\sum_b \xi_{bi}B_b)_{i\le c}$. For $L,r>0$ and
$(p,\xi)\in A\times Y$, let
\[
    \beta_L(p,\xi,r):=r^{-s}\mu\bigl(B_r(p)\cap
    \bigl\{z \ : \ |\Pi_\bV(z)-R^\xi(z)|>Lr^{m+\alpha}\bigr\}\bigr)
\]
and $\omega_{L,n}(p,\xi):=\sup\{\beta_L(p,\xi,r) \ : \ r\in\mathbb Q,\ 0<r<1/n\}$. The
integrand defining $\beta_L$ is jointly Borel in $(p,\xi,z)$, so
$\beta_L(\cdot,\cdot,r)$ and $\omega_{L,n}$ are Borel by Tonelli's theorem.
By~\cref{lem:rational_radii}, applied with $\varphi:=|\Pi_\bV-R^\xi|$ and
$\gamma:=m+\alpha$, we have
\begin{equation}\label{eq:rational_radii}
    \sup_{0<r<1/(2n)}\beta_{2^{m+\alpha}L}(p,\xi,r)
    \le2^s\omega_{L,n}(p,\xi).
\end{equation}
Hence $R^\xi$ satisfies \textup{(iii)} at $p$ if and only if
$\omega_{L,n}(p,\xi)\to0$ as $n\to\infty$ for some $L\in\mathbb N$ when
$\alpha>0$, and for every $L\in\{1/h:h\in\mathbb N\}$ when $\alpha=0$.
Condition \textup{(ii)} is closed in $(p,\xi)$, so the relation of pairs
$(p,\xi)\in A\times Y$ satisfying \textup{(i)--(iii)} is a Borel subset of
$\G\times Y$ with projection $A$. Applying~\cref{lem:borel_choice} with
$X=\G$ and $\mu=\cH^s\res A$ gives a Borel map
$p\mapsto \xi(p)$ such that
$R^{\xi(p)}$ satisfies \textup{(i)--(iii)} at almost every $p$. Deleting a null
set, we may assume this holds on $A$, and we write $R_p:=R^{\xi(p)}$ from now on
(by~\cref{rmk:jet_uniqueness}, this is the polynomial of the hypothesis).

\proofstep{3 (Compact pieces)} Set $g_{L,n}(p):=\omega_{L,n}(p,\xi(p))$,
Borel function on $A$. If $\alpha>0$, the Borel sets
$A^{(j)}:=\{p\in A \ : \ g_{j,n}(p)\to0\text{ as }n\to\infty\}$, $j\in\mathbb N$, cover
$A$ by Step~2; if $\alpha=0$, then $g_{1/h,n}\to0$ on $A$ for every
$h\in\mathbb N$. Apply~\cref{lem:compact_pieces} to $\mu$ on each $A^{(j)}$ with
the functions $g_{j,n}$ (respectively on $A$ with all $g_{1/h,n}$), with the
map $p\mapsto \xi(p)$, and with exponent $s$ and the original density measure
$\nu=\mu=\cH^s\res A$, unchanged when passing to $A^{(j)}$.
With~\eqref{eq:rational_radii}, this covers $A$ up to a null set by compact
sets $K$ on which:
\begin{enumerate}[label=(\alph*)]
\item $p\mapsto \xi(p)$ is continuous;
\item there are $c_1,r_1>0$ such that
      $\mu(B_r(p))\ge c_1r^s$ for every $p\in K$ and $0<r<r_1$;
\item if $\alpha>0$, there is $L>0$ such that
      \[
      \sup_{p\in K}r^{-s}\mu\bigl(B_r(p)\cap\{|\Pi_\bV-R_p|>Lr^{m+\alpha}\}\bigr)
      \to0
      \]
      as $r\downarrow0$. If $\alpha=0$, the same holds with $r^m$ in
      place of $r^{m+\alpha}$, for every $L>0$.
\end{enumerate}
For~\textup{(c)}, apply~\eqref{eq:rational_radii} with $\xi=\xi(p)$, taking
$L=j$ and absorbing $2^{m+\alpha}$ into the tube constant if $\alpha>0$, and
$L=1/h$ with $2^m/h\le\Lambda$ for a prescribed tube constant $\Lambda$ if
$\alpha=0$. Then, $g_{L,n}\to0$ uniformly on $K$.

\proofstep{4 (Whitney compatibility)} We compare $R_p$ and $R_q$
on a common set where both approximate $\Pi_\bV$. Positive lower density and
the small exceptional sets leave enough measure for~\cref{lem:normal_form_norming}
to turn this comparison into bounds for all derivatives.

Let $c_0,r_0$ be given
by~\cref{lem:normal_form_norming} with $K$ and $\theta=c_1/2$. For
distinct $p,q\in K$ set $\sigma:=d(p,q)$, and let $\eta_\sigma:=
(2^{m+\alpha}+1)L\sigma^\alpha$ if $\alpha>0$, with $L$ as in~\textup{(c)}, and
$\eta_\sigma:=\eta$ for a fixed $\eta>0$ if $\alpha=0$. Since
$B_\sigma(q)\subset B_{2\sigma}(p)$, condition~\textup{(c)}, applied at $p$ with
radius $2\sigma$ and at $q$ with radius $\sigma$ (for $\alpha=0$ with tube constants
$2^{-m-1}\eta$ and $\eta/2$), shows that for $\sigma$ smaller than some
$\sigma_*>0$ independent of $p,q$ (and depending on $\eta$ if $\alpha=0$) both exceptional sets have measure at most
$(c_1/4)\sigma^s$. By~\textup{(b)}, their complement in $A\cap B_\sigma(q)$ is a
Borel set $Z_{p,q}$ with $\mu(Z_{p,q})\ge(c_1/2)\sigma^s$, and
\begin{equation}\label{eq:both_approximations_good}
    |\Pi_{\bV}(z)-R_p(z)|+|\Pi_{\bV}(z)-R_q(z)| \le \eta_\sigma\sigma^m,
    \qquad z\in Z_{p,q}.
\end{equation}
Applying~\cref{lem:normal_form_norming} at $x=q$, with $r=\sigma$ and
$Z=Z_{p,q}$, to $R_q-R_p\in[\mathcal N_{\bV}^m]^c$ gives, for
$\sigma<\delta_0:=\min\{\sigma_*,r_0,r_1\}$, the following derivative inequality:
\begin{equation}\label{eq:trace_selection_whitney_estimate_alpha}
    |\mbX^J(R_q-R_p)(q)| \le c_0^{-1}\eta_\sigma\,\sigma^{m-d(J)},
    \qquad d(J)\le m.
\end{equation}

\proofstep{5 (Whitney conditions)} Since $R_p$ has homogeneous degree at most $m$, it
is its own Taylor polynomial at $p$. Hence the polynomial determined by the
jet $\{A_I(p)\}_{d(I)\le m}$ is $R_p$, and it satisfies
\[
    A_J(q)-\mbX^JP_m\bigl[\{A_I\}_{d(I)\le m}\bigr](p,q)=\mbX^J(R_q-R_p)(q).
\]
By~\textup{(a)}, the functions $A_J$ are continuous, hence bounded, on $K$. If
$\alpha=0$ and $\eps>0$ is prescribed, take $\eta:=c_0\eps$,
then~\eqref{eq:trace_selection_whitney_estimate_alpha} is the little-$o$
condition~\eqref{eq:whitney_little_o_condition} for all $p\ne q\in K$ with
$d(p,q)<\delta_0$, hence on $K\cap B_{\delta_0/2}(\bar x)$ for every
$\bar x\in K$, uniformly in $\bar x$. If
$\alpha>0$, it gives the bound of~\cref{def:lip_jet} with $\gamma=m+\alpha$
for nearby pairs, and for the remaining pairs the left-hand side is bounded.
Enlarging the constant gives the jet assertion.

\proofstep{6 (Conclusion)} Let $K$ be one of these compact sets.
\Cref{thm:whitney_extension} if
$\alpha=0$, and~\cref{thm:whitney_alpha_holder} if $\alpha>0$, applied
componentwise, give $\Xi\in[\mathbf C^{m,\alpha}_{\G}(\G)]^c$ with
$\mbX^J\Xi=A_J$ on $K$ for $d(J)\le m$. The polynomial map $\Pi_\bV$ is smooth,
so $F:=\Pi_{\bV}-\Xi\in[\mathbf C^{m,\alpha}_{\G}(\G)]^c$. For $p\in K$ we have
$F(p)=\Pi_\bV(p)-R_p(p)=0$. Moreover, each $V_i$ is a linear
combination of $X_1,\ldots,X_\ell$ and $m\ge1$, so
$V_i\Xi(p)=V_iR_p(p)=0$, because $V_iR_p\equiv0$; since $V_i\Pi_\bV=e_i$, we
get $V_iF(p)=e_i$. By~\cref{lem:level_set}, $K$ lies in a $k$-dimensional
$\G$-regular surface of class $(m,\alpha)$, and countably many such sets $K$
cover $A$ up to an $\cH^s$-null set.
\end{proof}

\begin{remark}\label{rmk:jet_uniqueness}
At every $p\in A$ with $\Theta_*^s(A,p)>0$ there is at most one polynomial
satisfying \textup{(i)--(iii)} of~\cref{lem:trace_normalized_selection}. In
particular, the selected jet does not depend on the surfaces from which the
polynomials are obtained in~\cref{thm:higher_order}. Indeed, let $R,R'$
satisfy \textup{(i)--(iii)} at such $p$, put $P:=R-R'\in[\mathcal N_{\bV}^m]^c$,
and choose $\theta>0$ with $\cH^s(A\cap B_r(p))\ge3\theta r^s$ for small $r$.
If $\alpha>0$, set $\eta_r:=(L+L')r^\alpha$, where $L,L'$ are tube constants
for $R,R'$; if $\alpha=0$, fix $\eta>0$, set $\eta_r:=\eta$, and use
\textup{(iii)} with tube constant $\eta/2$ for both. For small $r$, removing
the two exceptional sets from $A\cap B_r(p)$ leaves a set containing a Borel
set $Z_r$ with $\cH^s(Z_r)\ge\theta r^s$ and $|P|\le\eta_rr^m$ on $Z_r$.
\Cref{lem:normal_form_norming} with $K=\{p\}$ gives
$r^{d(J)}|\mbX^JP(p)|\le c_0^{-1}\eta_rr^m$ for $d(J)\le m$. Letting
$r\downarrow0$, and then $\eta\downarrow0$ if $\alpha=0$, shows that
$\mbX^JP(p)=0$ for $d(J)\le m$, so $P\equiv0$ by the uniqueness
in~\eqref{eq:Taylor_cond}.
\end{remark}

%%%%%%%%%%%%%%%%%%%%%%%%%%%%%%%%%%%%%%%%%%%%%%%%%%%%%%%%%%%%%%%%%%%%%%%%%%%%%%%%%%%%%%%%%%%%%%%%%%%%%%%
\section{Higher-order rectifiability}\label{sec:higher-order-rectifiability}
%%%%%%%%%%%%%%%%%%%%%%%%%%%%%%%%%%%%%%%%%%%%%%%%%%%%%%%%%%%%%%%%%%%%%%%%%%%%%%%%%%%%%%%%%%%%%%%%%%%%%%%

We now prove the order-one criterion and the higher-order one (\cref{thm:higher_order}), and combine
them to obtain~\cref{thm:main}. The
higher-order information is encoded by normal-form ambient jets for the split
projection $\Pi_{\bV}$ on first-order graph patches, and the rectifying
surfaces are non-critical level sets of ambient $\mathbf C^{m,\alpha}_{\G}$
functions.

\subsection{Order-one rectifiability}\label{sec:1-rect}

The following criterion is the general-Carnot analogue
of~\cite{MSSC10}*{Theorem~3.15}.
Antonelli and Merlo~\cite{AntonelliMerlo2022Representation} characterize
rectifiable measures in Carnot groups through tangent measures that are
multiples of Haar measures of homogeneous subgroups, and leave approximate
tangents in the sense of~\cite{mattila}*{Definition~15.17} to future work.
For $\cH^{k_{\m}}\res E$ with positive lower
density, tangent measures that are multiples of the Haar measure of one
$\bW_p\in\cV(\G,k)$ imply the tangent condition of~\cref{thm:c1_rec}, which
constrains only the directions in which $E$ concentrates at small scales.
Marstrand--Mattila-type criteria, in which flat tangent measures may vary with
the scale, are proved in~\cite{Merlo2023Marstrand} in codimension one and
in~\cite{AntonelliMerlo2022Marstrand} for tangent subgroups with a normal
complement. Proper vertical subgroups containing the center, including all proper vertical
subgroups of Heisenberg groups, have no normal complement, since every nontrivial normal
subgroup of a nilpotent group meets the center.
The proof of the criterion uses local horizontal complements, measurable tangent
selection, compact-piece uniformization, and Whitney extension.
First, we state a technical lemma concerning the local stability of the complements.

\begin{lemma}\label{lem:local_complements}
Let $\bW_0\in\cV(\G,k)$ and let $\bV_0$ be a horizontal complement. Then $\G=\bW\rtimes\bV_0$ for every $\bW$ in a neighborhood of
$\bW_0$ in $\cV(\G,k)$. Additionally, under the identification
$\bW\leftrightarrow\bW^1$, the space $\cV(\G,k)$ is an open subset of the
Grassmannian of $(k-N+\ell)$-dimensional subspaces of $\cg_1$. In particular,
it is Polish.
\end{lemma}

\begin{proof}
As in the proof of~\cref{prop:vertical_splitting}, write
$\cg_1=\bW_0^1\oplus\Lie(\bV_0)$ with $\Lie(\bV_0)$ abelian.
Nearby $(k-N+\ell)$-dimensional subspaces are graphs of linear maps
$L:\bW_0^1\to\Lie(\bV_0)$, hence remain complementary to $\Lie(\bV_0)$.
Thus transversality is open, and for every such
transverse $S$~\cref{prop:vertical_splitting} gives $\G=\bW_S\rtimes\bV_0$, where
$\bW_S:=\exp(S\oplus\cg_2\oplus\cdots\oplus\cg_\iota)$ is vertical and
$k$-dimensional. Since every element of $\cV(\G,k)$ equals $\bW_S$ with
$S=\bW^1$, both assertions follow.
\end{proof}

Fix $1\le k<N$. By~\cref{lem:local_complements}, $\cV(\G,k)$ is an open subset of a compact
Grassmannian, hence locally compact and second countable, so it is covered by
countably many open sets $\mathcal U_a$, indexed by $a\in\mathbb N$, with
compact closures in $\cV(\G,k)$, each with a horizontal subgroup $\bV_a$ such that
$\G=\bW\rtimes\bV_a$ for every $\bW\in\overline{\mathcal U_a}$. Identify $\bV_a$ with $\R^{N-k}$ through a fixed basis and, for
$\bW\in\overline{\mathcal U_a}$, let $\varpi_{a,\bW}:\G\to\R^{N-k}$ be the
coordinate representation of the split projection onto $\bV_a$.
By~\cref{thm:semidirect_fact}, $\varpi_{a,\bW}$ is an h-homomorphism, and
\begin{equation}\label{eq:c1_projection_kernel}
    \ker\varpi_{a,\bW}=\bW,
    \qquad
    \varpi_{a,\bW}|_{\bV_a}=\operatorname{Id}_{\bV_a}.
\end{equation}
The map $\varpi_{a,\bW}$ depends only on the first-layer coordinates, on
which it is the linear projection onto $\Lie(\bV_a)$ along $\bW^1$. Hence the
projections of the splittings $\bW\rtimes\bV_a$ depend continuously on
$(\bW,x)$, and compactness of $\overline{\mathcal U_a}$ and homogeneity give
constants $0<c_a\le1 \le C^*_a$, and a nondecreasing function
$\omega_a$ with $\omega_a(t)\to0$ as $t\downarrow0$, such that, for all
$x\in\G$, $v\in\bV_a$, and $\bW,\bW'\in\overline{\mathcal U_a}$,
\begin{equation}\label{eq:c1_projection_lipschitz}
 \begin{gathered}
    C_{\bW,\bV_a}\le C^*_a,\qquad c_a|v|\le\|v\|\le c_a^{-1}|v|,\\
    |\varpi_{a,\bW}(x)-\varpi_{a,\bW'}(x)|\le\omega_a(\rho(\bW,\bW'))\lVert x\rVert,
 \end{gathered}
\end{equation}
where $C_{\bW,\bV_a}$ is the constant of~\cref{thm:semidirect_fact} for the
splitting $\bW\rtimes\bV_a$, which depends continuously on $\bW$, and $|v|$
is the Euclidean norm of the coordinates. The comparison holds because
$\bV_a$ is horizontal, on which $\|\cdot\|$ is homogeneous of
degree one.

\begin{theorem}\label{thm:c1_rec}
Let $1\le k\le N$, and let $E\subset\G$ be $\cH^{k_{\m}}$-measurable with $\cH^{k_{\m}}(E)<\infty$ and $\Theta_*^{k_{\m}}(E,p)>0$ for $\cH^{k_{\m}}$-almost every $p\in E$. Suppose that for $\cH^{k_{\m}}$-a.e. $p\in E$ there exists $\bW_p\in\cV(\G,k)$ such that
\begin{equation}\label{eq:c1_approximate_tangent_condition}
    \Tan^{*k_{\m}}(\cH^{k_{\m}}\res E,p)\subseteq\bW_p.
\end{equation}
Then $E$ is $(k,\G)$-rectifiable of class $(1,0)$ in the sense of~\cref{def:rec}.
\end{theorem}

\begin{proof}
If $k=N$, the conclusion is trivial because $\G$ itself is an $N$-dimensional regular surface. Hence assume that $k<N$, and set $s:=k_{\m}$, $c:=N-k$, and $\mu:=\cH^s\res E$. We may also assume that $\mu(E)>0$ and, replacing $E$ by a Borel representative (which exists because $\cH^s$ is Borel regular and $\cH^s(E)<\infty$), that $E$ is Borel. Then, $\mu$ is a finite Radon measure on the Polish space $\G$ by~\cite{Kechris1995}*{Theorems~17.10--17.11}.

\proofstep{1 (Borel choice of approximate tangents)} We use the charts fixed in the discussion above. For
$j,n\in\mathbb N$, $p\in\G$, and $\bW\in\mathcal U_a$, set
\[
\zeta_{a,j,n}(p,\bW):=\sup\left\{r^{-s}\mu\left(B_r(p)\cap
\left\{q \ : \ |\varpi_{a,\bW}(p^{-1}q)|> \tfrac{r}{j}\right\}\right) \ : \ 
r\in\mathbb Q,\ 0<r<\tfrac{1}{n}\right\}.
\]
For fixed rational $r$, the quantity inside the supremum is $r^{-s}$ times
the $\mu$-integral in $q$ of the jointly Borel function
$\mathbf 1_{\{d(p,q)<r\}}\mathbf 1_{\{|\varpi_{a,\bW}(p^{-1}q)|>r/j\}}$, so
every $\zeta_{a,j,n}$ is Borel by Tonelli's theorem. If $\eta\ge2/j$ and
$0<r<1/(2n)$, then~\cref{lem:rational_radii}, applied with
$\varphi:=|\varpi_{a,\bW}(p^{-1}\,\cdot)|$, $\gamma:=1$, and $L:=1/j$, gives
\[
 r^{-s}\mu\bigl(B_r(p)\cap\bigl\{q \ : \ |\varpi_{a,\bW}(p^{-1}q)|>\eta r\bigr\}\bigr)
 \le 2^s\zeta_{a,j,n}(p,\bW).
\]
Hence $\zeta_{a,j,n}(p,\bW)\to0$ as $n\to\infty$ for every $j$ if and only if
\begin{equation}\label{eq:c1_concentration}
    \lim_{r\downarrow0}r^{-s}\mu\bigl(B_r(p)\cap
    \bigl\{q \ : \ |\varpi_{a,\bW}(p^{-1}q)|>\eta r \bigr\}\bigr)=0
    \qquad\text{for every }\eta>0.
\end{equation}
On $\bV_a$ the homogeneous norm is comparable to the Euclidean norm of the
coordinates, so by the implication \textup{(i)$\Rightarrow$(iii)}
of~\cref{lem:tangent_condition}, applied to the splitting $\bW_p\rtimes\bV_a$
(not to a complement of $\bW_p$ fixed in advance), \eqref{eq:c1_concentration}
holds with $\bW=\bW_p$ whenever $\bW_p\in\mathcal U_a$
satisfies~\eqref{eq:c1_approximate_tangent_condition}. Therefore the Borel set
\[
 \mathcal R:=\bigcap_{j,i\ge1} \bigcup_{n_0\ge1} \bigcap_{n\ge n_0}
 \bigl\{(p,a,\bW) \ : \ \zeta_{a,j,n}(p,\bW)< \tfrac1i\bigr\}
 \subset\G\times\mathcal Y, \quad \mathcal Y:=\bigsqcup_a\bigl(\{a\}\times\mathcal U_a\bigr),
\]
has projection of full $\mu$-measure, and~\cref{lem:borel_choice} gives a
Borel map $p\mapsto(a(p),\bW_p)$ with $(p,a(p),\bW_p)\in\mathcal R$ for
$\mu$-a.e.\ $p$. Thus $a(p)$ selects a chart containing the chosen $\bW_p$ and
it becomes constant on each compact piece below. To simplify the notation, we use this choice of $\bW_p$ from now on and write
$\varpi_p:=\varpi_{a(p),\bW_p}$, so that~\eqref{eq:c1_concentration} holds
with $\varpi_p$ in place of $\varpi_{a,\bW}$ at $\mu$-a.e.\ $p$. By the
implication \textup{(iii)$\Rightarrow$(i)} of~\cref{lem:tangent_condition},
applied to $\bW_p\rtimes\bV_{a(p)}$, the chosen $\bW_p$
satisfy~\eqref{eq:c1_approximate_tangent_condition} at $\mu$-a.e.\ $p$.

\proofstep{2 (Compact pieces and uniform flatness)} Remove the null set
where the lower-density hypothesis or~\eqref{eq:c1_concentration} fails.
Apply~\cref{lem:compact_pieces} on each chart to the original measure
$\mu=\cH^s\res E$, the map $p\mapsto\bW_p$, the functions
$p\mapsto\zeta_{a,j,n}(p,\bW_p)$ for every $j$, and the lower-density bound
with $\nu=\mu$ and exponent $s$. Together with the rational-radius comparison
in \textbf{Step~1}, this gives countably many compact sets $K$ covering $E$ up to a
null set. On each $K$, the chart index $a$ is constant, $p\mapsto\bW_p$ is
continuous, and there are a constant $b>0$ and a radius $r_0>0$ such that
\[
 \mu(B_t(q))\ge bt^s\qquad \text{for all } q\in K \ \text{ and } \ 0<t<r_0.
\]
Moreover, for every $\eta>0$,
\[
 \sup_{p\in K}r^{-s}\mu\bigl(B_r(p)\cap
 \{z \ : \ |\varpi_p(p^{-1}z)|>\eta r\}\bigr)\longrightarrow0.
\]
Fix such a $K$, let $\bV:=\bV_a$, and write $L:=c_a^{-1}C^*_a$. Then
$|\varpi_p(g)|\le L\|g\|$ uniformly for $p\in K$.
Fix $\eps>0$ and set $\delta:=\min\{1/2,\eps/(4L)\}$.
If $p,q\in K$, $r:=d(p,q)>0$, and
$|\varpi_p(p^{-1}q)|>\eps r$, the homomorphism property gives
the inclusion
\[
 B_{\delta r}(q)\subset B_{2r}(p)\cap
 \bigl\{z \ : \ |\varpi_p(p^{-1}z)|>(\eps/4)(2r)\bigr\}.
\]
For small $r$, the left-hand side has $\mu$-measure at least
$b\delta^sr^s$, whereas the right-hand side has measure $o(r^s)$ uniformly
in $p\in K$. This is impossible for sufficiently small $r$, uniformly in
$p,q$, and therefore
\begin{equation}\label{eq:uniform_flatness}
 \lim_{r\downarrow0}\sup\left\{
 \frac{|\varpi_p(p^{-1}q)|}{d(p,q)} \ : \ p,q\in K,\ 0<d(p,q)<r\right\}=0.
\end{equation}
Finally, $p\mapsto\varpi_p$ is uniformly
continuous on $K$ by~\eqref{eq:c1_projection_lipschitz}.

\proofstep{3 (Whitney extension)} If $X_1,\ldots,X_\ell$ is the fixed horizontal
basis (i.e., a basis of the horizontal layer), prescribe the $\R^c$-valued first-order jet on $K$ by
\[
    F_0(p):=0,
    \qquad
    F_{J_i}(p):=\varpi_p(\exp X_i),\quad i=1,\ldots,\ell,
\]
where $J_i$ is the multi-index with $\mbX^{J_i}=X_i$. The jet is bounded on $K$, because $p\mapsto\varpi_p$ is continuous and $K$ is compact. The first-order polynomial of this jet at $p$ (in other words, the polynomial $P_1$
of~\cref{subsec:whitney}) is
\[
    P_p(q)=\varpi_p(p^{-1}q),
\]
and, for every horizontal basis vector $X_i$,
\[
    X_iP_p(q)=\varpi_p(\exp X_i).
\]
Consequently, for $p,q\in K$,
\[
    \bigl|F_{J_i}(q)-X_iP_p(q)\bigr|
    =\bigl|\varpi_q(\exp X_i)-\varpi_p(\exp X_i)\bigr|=o(1)
\]
uniformly as $d(p,q)\to0$, by the uniform continuity of $p\mapsto\varpi_p$. This is the derivative part of the Whitney condition. For the value part,~\eqref{eq:uniform_flatness} gives
\[
    |F_0(q)-P_p(q)|=|\varpi_p(p^{-1}q)|=o(d(p,q)).
\]
Thus the jet satisfies the Whitney compatibility conditions of order one.
Applying~\cref{thm:whitney_extension} componentwise gives a map
\[
    F\in[\mathbf C^1_\G(\G)]^c
\]
such that $F=0$ and $d_HF(p)=\varpi_p$ on $K$. For every $p\in K$, we have $F(p)=0$, and~\eqref{eq:c1_projection_kernel}
gives $d_HF(p)|_{\bV}=\operatorname{Id}_{\bV}$, that is, $V_iF(p)=e_i$ in the
fixed coordinates on $\bV$. By~\cref{lem:level_set}, $K$ lies in a
$k$-dimensional $\G$-regular surface of class $1$. Taking the countable
union over the sets $K$ concludes the proof.
\end{proof}

\subsection{Proof of the main theorems}\label{sec:higher-order-rect}

The proof of~\cref{thm:main} obtains approximate tangents from the tubular
approximation at single points.

\begin{lemma}\label{lem:tube_to_tangent}
Let $1\le k\le N$, $m\in\N$, $\alpha\in[0,1]$, and $s:=k_{\m}$. Let
$E\subset\G$ be $\cH^s$-measurable and locally $\cH^s$-finite, let $p\in E$,
and let $S$ be a $k$-dimensional $\G$-regular surface of class $(m,\alpha)$
through $p$ such that~\eqref{eq:main_assumption} holds with $S_p=S$ for some
$\lambda\ge0$ if $\alpha>0$, and~\eqref{eq:main_assumption_alpha_zero} holds
with $S_p=S$ for every $\lambda>0$ if $\alpha=0$. Let $\bV_p$ be a horizontal
complement of $\bW_p:=T_pS\in\cV(\G,k)$. Then
\begin{equation}\label{eq:tube_to_tangent}
\lim_{r\downarrow0}r^{-s}\cH^s\Bigl(E\cap B_r(p)\cap
\bigl\{q \ : \ \lVert\pi_{\bV_p}(p^{-1}q)\rVert>\eps r\bigr\}\Bigr)=0
\qquad\text{for every }\eps>0,
\end{equation}
and consequently $\Tan^{*s}(\cH^s\res E,p)\subseteq\bW_p$.
\end{lemma}

\begin{proof}
We first show the ordinary first-order tangency of $S$ to $p\cdot\bW_p$. If $F$ is a local defining function for $S$, then for $b\in S$ close to $p$ the first-order Taylor formula gives
\[
    d_HF(p)(p^{-1}b)=o(d(p,b)).
\]
The h-homomorphism $d_HF(p)$ has kernel $\bW_p$ and restricts to an isomorphism
\[
    L_p:=d_HF(p)|_{\bV_p}:\bV_p\longrightarrow\R^{N-k}.
\]
Since $d_HF(p)=L_p\circ\pi_{\bV_p}$, the boundedness of $L_p^{-1}$, together with the comparability of the homogeneous and Euclidean norms on the horizontal subgroup $\bV_p$, implies
\begin{equation}\label{eq:approximating_surface_linear_tangent}
    \lVert\pi_{\bV_p}(p^{-1}b)\rVert=o(d(p,b))
    \qquad\text{as }b\to p\text{ in }S.
\end{equation}
Fix $\eps>0$ and $0<\delta<\eps/(4C_{\bW_p,\bV_p})$. If $q\in E\cap B_r(p)$ and $d(q,S)\le\delta r$, choose $b\in S$ with $d(q,b)\le2\delta r$. Then $d(p,b)\le(1+2\delta)r$, so~\eqref{eq:approximating_surface_linear_tangent} is uniform for all such choices of $b$. Since $\pi_{\bV_p}$ is an h-homomorphism,~\eqref{eq:approximating_surface_linear_tangent} and~\eqref{eq:piVcontin} give, for all sufficiently small $r$,
\[
    \lVert\pi_{\bV_p}(p^{-1}q)\rVert
    \le \lVert\pi_{\bV_p}(p^{-1}b)\rVert
       +\lVert\pi_{\bV_p}(b^{-1}q)\rVert
    \le o(r)+2C_{\bW_p,\bV_p}\delta r<\eps r.
\]
Consequently, for all sufficiently small $r$,
\[
E\cap B_r(p)\cap
\bigl\{q \ : \ \lVert\pi_{\bV_p}(p^{-1}q)\rVert>\eps r\bigr\} \subseteq E\cap B_r(p)\cap\{q \ : \ d(q,S)>\delta r\}.
\]
The set on the right has vanishing normalized $\cH^s$-measure. Indeed, if $\alpha>0$, then $\lambda r^{m+\alpha}\le\delta r$ for all small $r$; if $\alpha=0$ and $m=1$, use~\eqref{eq:main_assumption_alpha_zero} with tube constant $\delta$; and if $\alpha=0$ and $m>1$, use any fixed tube constant and $r^m=o(r)$. This proves~\eqref{eq:tube_to_tangent}, and~\cref{lem:tangent_condition} gives the inclusion of tangent cones.
\end{proof}

The proof of the higher-order criterion,~\cref{thm:higher_order}, is divided into four steps. Step~1
refines the first-order cover to graph patches. Step~2
restricts to the density part $A$ of one patch (here positive lower density on each patch is obtained, not assumed). Step~3 transfers the tube
estimate to normal-form ambient jets for $\Pi_{\bV}$ on $A$. Step~4 selects
compatible jets and constructs the regular surfaces as non-critical zero sets.

\begin{proof}[Proof of~\cref{thm:higher_order}]
If $k=N$, the conclusion is trivial, since $\G$ is an $N$-dimensional regular surface. Let $k<N$. Replacing $E$ by a Borel representative changes no approximation or rectifiability conclusion, so we may assume that $E$ is Borel.

\proofstep{1} By hypothesis, $E$ is $(k,\G)$-rectifiable of class $(1,0)$.
Refining the resulting regular surfaces by~\cref{thm:intrinsic_surface_c1}
and second countability, we obtain a countable cover, up to an
$\cH^{k_{\m}}$-null set, by split-regular intrinsic $\mathbf C_H^1$ graph
patches. Choose these as relatively compact subpatches of larger regular
graph charts. The local upper Ahlfors bound
in~\cite{AGV2022}*{Proposition~2.17}, applied on their compact closures in
the larger charts, makes each patch have finite $\cH^{k_{\m}}$-measure.

\proofstep{2} Fix one split-regular graph patch
\[
    M=\Phi_\phi(U), \qquad \Phi_\phi(w)=w\cdot \phi(w), \qquad \phi:U\subset\bW\to\bV,
\]
from this cover, where $\G=\bW\rtimes\bV$ and $\bV$ is the horizontal complement selected for the patch. We are thus in the setting of~\cref{subsec:graph_jets}. It is enough to prove that $E\cap M$ is covered, up to an $\cH^{k_{\m}}$-null set, by $\G$-regular surfaces of class $(m,\alpha)$. Set $c:=N-k$, $s:=k_{\m}$, $\nu:=\cH^s\res M$, and
\[
    A:=\left\{a\in E\cap M \ : \ \Theta^{*s}\bigl(\cH^{s}\res(M\setminus E),a\bigr)=0 \right\}.
\]
The measure $\nu$ is locally doubling
and locally Ahlfors $s$-regular by~\cite{AGV2022}*{Proposition~2.17}. The
differentiation theorem~\cite{AGV2022}*{Proposition~2.3(ii)--(iii)}, applied
to $\nu\res(M\setminus E)$ relative to $\nu$, therefore gives
\[
   \frac{\nu((M\setminus E)\cap B_r(a))}{\nu(B_r(a))}\longrightarrow0, \qquad \text{for $\nu$-almost every } a\in E\cap M,
\]
and the local upper Ahlfors bound
implies $r^{-s}\nu((M\setminus E)\cap B_r(a))\to0$. Thus $A$ has full
$\cH^s$-measure in $E\cap M$, and it is Borel because the defining limsup can
be taken over rational radii. It is enough to cover $A$. Since
$M\setminus A\subseteq(M\setminus E)\cup((E\cap M)\setminus A)$, where the
second set is $\cH^s$-null and the first has zero upper $s$-density at every
point of $A$,
\begin{equation}\label{eq:density_part_full_in_graph}
    \lim_{r\downarrow0}r^{-s}\cH^{s}\bigl((M\setminus A)\cap B_r(a)\bigr)=0
    \qquad\text{for every }a\in A.
\end{equation}
Fix $a\in A$ and consider the rescalings
\[
    M_{a,r}:=\delta_{1/r}(a^{-1}M),
    \qquad A_{a,r}:=\delta_{1/r}(a^{-1}A).
\]
\Cref{lem:graph_blowup} gives
$\cH^s\res M_{a,r}\overset{*}{\rightharpoonup}\cH^s\res T_aM$. Moreover, for every $R<\infty$,~\eqref{eq:density_part_full_in_graph} gives
\[
 \cH^s\bigl((M_{a,r}\setminus A_{a,r})\cap B_R(e)\bigr)
 =r^{-s}\cH^s\bigl((M\setminus A)\cap B_{Rr}(a)\bigr)\longrightarrow0.
\]
Consequently
\begin{equation}\label{eq:blowup_of_A}
    \cH^s\res A_{a,r}\overset{*}{\rightharpoonup}\cH^s\res T_aM.
\end{equation}
Homogeneous subgroup spheres have zero Haar measure
by~\cite{AGV2022}*{Lemma~3.5}. Hence the continuity-set conclusion
of~\cite{AFP}*{Proposition~1.62(b)}, applied to $B_1(e)$, shows that
$r^{-s}\cH^s(A\cap B_r(a))$ converges to the positive finite number
$\cH^s(T_aM\cap B_1(e))$. In particular, $A$ has positive lower $s$-density
at every point.

\proofstep{3} We show that a.e. point of $A$ gives a polynomial $R_a$ satisfying conditions~\textup{(i)--(iii)} of~\cref{lem:trace_normalized_selection}, with $\Pi_\bV$ as in~\cref{subsec:graph_jets}. Fix $a\in A$ at which the defining $\mathfrak S_{m,\alpha}(k)$ approximation property holds, and let $S_a$ be an approximating $k$-dimensional $\G$-regular surface of class $(m,\alpha)$ through $a$. We first identify the tangent of $S_a$.

Let $\bV'$ be a horizontal complement of $T_aS_a$. \Cref{lem:tube_to_tangent}, applied to
$S_a$ and to $A$, whose exceptional sets are contained in those of $E$, gives
\[
    \lim_{r\downarrow0}r^{-s}\cH^s\Bigl(A\cap B_r(a)\cap
    \bigl\{q \ : \ \lVert\pi_{\bV'}(a^{-1}q)\rVert>\eps r\bigr\}\Bigr)=0
    \qquad\text{for every }\eps>0.
\]
By the homogeneity of $\pi_{\bV'}$, the quantity under the limit equals
$\cH^s(A_{a,r}\cap O_\eps)$, where
$O_\eps:=B_1(e)\cap\{z:\lVert\pi_{\bV'}(z)\rVert>\eps\}$ is open.
By~\eqref{eq:blowup_of_A} and the open-set bound
of~\cite{AFP}*{Proposition~1.62(a)}, $\cH^s(T_aM\cap O_\eps)=0$. Since
$\cH^s\res T_aM$ is a Haar measure of $T_aM$, it charges every nonempty
relatively open subset of $T_aM$, so $T_aM\cap O_\eps=\emptyset$ for every
$\eps>0$. Thus $T_aM\cap B_1(e)\subseteq\ker\pi_{\bV'}=T_aS_a$, and
$T_aM\subseteq T_aS_a$ by dilation invariance. Both are vertical homogeneous
subgroups of topological dimension $k$, so necessarily
\begin{equation}\label{eq:tangent_B_equals_M}
    T_aS_a=T_aM.
\end{equation}
Because the patch $M$ is split-regular over the fixed splitting $\G=\bW\rtimes\bV$, its tangent $T_aM$ is transverse to $\bV$. Hence~\eqref{eq:tangent_B_equals_M} gives
\[
    T_aS_a \cap \bV=\{0\}.
\]
Since $T_aS_a$ is vertical and the dimensions add up, this is the splitting
\[
    \G=T_aS_a\rtimes\bV=T_aM\rtimes\bV.
\]
%We use only an ambient defining function of $S_a$, not a graph representation of it. 
Now, by~\cref{lem:surface_normal_form_jet} applied to $S_a$ at $a$, we get polynomials $R_a\in[\mathcal N_\bV^m]^c$ with $R_a(a)=\Pi_{\bV}(a)$
that satisfy, for some $r_a>0$,
\[
    |\Pi_{\bV}(b)-R_a(b)| \le C_a d(a,b)^{m+\alpha}
    \qquad\text{for } b\in S_a\cap B_{r_a}(a)
\]
when $\alpha>0$. When $\alpha=0$, the right-hand side is
$\widetilde\omega_a(d(a,b))d(a,b)^m$, where
$\widetilde\omega_a$ may be chosen nondecreasing and tends to zero at zero.

We claim that $R_a$ also satisfies condition~\textup{(iii)}
of~\cref{lem:trace_normalized_selection} at $a$. Set $D_a:=\Pi_{\bV}-R_a$. Smooth ambient
polynomials are locally Lipschitz with respect to the homogeneous distance. Indeed,
on a compact set $K$, one has $|x-y|\le C_Kd(x,y)$, since
$y-x=h+P(x,h)$ for $h:=x^{-1}y$, with $P$ polynomial and $P(x,0)=0$
by~\eqref{eq:BCH}, and $|h|\le C\|h\|$ for $\|h\|\le1$. Hence, after
decreasing $r_a$, there is $L_a<\infty$ such that $D_a$ is
$L_a$-Lipschitz on $B_{r_a}(a)$. Fix $0<2r<r_a$ and
$q\in A\cap B_r(a)$. Choose $b_j\in S_a$ such that
\[
    d(q,b_j)\longrightarrow d(q,S_a).
\]
Since $a\in S_a$, one has $d(q,S_a)\le d(q,a)<r$. Thus, for all sufficiently
large $j$, $d(q,b_j)<r$ and $d(a,b_j)<2r<r_a$. For $\alpha>0$,
\[
    |D_a(q)|\le L_a d(q,b_j)+C_a(2r)^{m+\alpha}.
\]
Letting $j\to\infty$ and enlarging a constant $C\ge1$ depending on $a$
gives
\[
    |\Pi_{\bV}(q)-R_a(q)|
    \le C d(q,S_a)+C r^{m+\alpha}
\]
when $\alpha>0$. For $\alpha=0$, the same argument gives
\[
    |\Pi_{\bV}(q)-R_a(q)|
    \le C d(q,S_a)+\omega_a(r)r^m,
    \qquad \omega_a(r):=2^m\widetilde\omega_a(2r)\longrightarrow0.
\]
This also covers $d(q,S_a)=0$, without requiring a nearest point in $S_a$.
Since $A\subset E$ and $S_a$ approximates $E$ at $a$ in the sense of~\cref{def:special_class}, this gives condition~\textup{(iii)}. Indeed, suppose first that $\alpha>0$, and let $\lambda(a)$ be the tube constant supplied by~\cref{def:special_class}. Choose $\Lambda>C(\lambda(a)+1)$. The preceding estimate then gives, for all sufficiently small $r$,
\[
\begin{split}
&A\cap B_r(a)\cap
\bigl\{q:|\Pi_{\bV}(q)-R_a(q)|>\Lambda r^{m+\alpha}\bigr\}\\
&\hspace{5em}\subseteq E\cap B_r(a)\cap
\bigl\{q:d(q,S_a)>\lambda(a)r^{m+\alpha}\bigr\}.
\end{split}
\]
Consequently,
\[
    \lim_{r\downarrow0}r^{-k_{\m}}\cH^{k_{\m}}\Bigl( A\cap B_r(a)\cap \{q \ : \ |\Pi_{\bV}(q)-R_a(q)|>\Lambda r^{m+\alpha}\} \Bigr)=0.
\]
When $\alpha=0$, fix any $\Lambda>0$, choose the tube constant $\lambda<\Lambda/(2C)$, and then take $r$ so small that $\omega_a(r)<\Lambda/2$. The same argument gives
\[
\begin{split}
&A\cap B_r(a)\cap
\bigl\{q:|\Pi_{\bV}(q)-R_a(q)|>\Lambda r^m\bigr\}\\
&\hspace{7em}\subseteq E\cap B_r(a)\cap
\bigl\{q \ : \ d(q,S_a)>\lambda r^m\bigr\},
\end{split}
\]
and the normalized measure of the latter set tends to zero. Since $\Lambda$ was arbitrary, condition~\textup{(iii)} holds also when $\alpha=0$. Thus $R_a$ satisfies \textup{(i)--(iii)} at $\cH^{k_{\m}}$-a.e. $a\in A$.

\proofstep{4} The set $A$ is Borel by Step~2 and has finite $\cH^{k_{\m}}$-measure because $A\subset E$. By Steps~2 and~3, $A$ satisfies the hypotheses of~\cref{lem:trace_normalized_selection}, which yields countably many $k$-dimensional $\G$-regular surfaces $\{S_j\}_{j=1}^\infty$ of class $(m,\alpha)$ such that
\[
    \cH^{k_{\m}}\left( A\setminus\bigcup_{j=1}^\infty S_j \right)=0.
\]
Because $E\cap M$ differs from $A$ by an $\cH^{k_{\m}}$-null set, the same family covers $E\cap M$ up to null measure. Summing over the countably many first-order graph patches $M$ obtained in Step~1 gives a countable family of $k$-dimensional $\G$-regular surfaces of class $(m,\alpha)$ covering $E$ up to an $\cH^{k_{\m}}$-null set, concluding the proof.
\end{proof}

As a consequence, we can finally obtain the proof of our main result.

\begin{proof}[Proof of~\cref{thm:main}]
For $\cH^{k_{\m}}$-almost every $p\in E$, let $S_p$ be the $\G$-regular surface of class $(m,\alpha)$ through $p$ given by~\cref{def:special_class}, and set $\bW_p:=T_pS_p\in\cV(\G,k)$. By~\cref{lem:tube_to_tangent}, $\Tan^{*k_{\m}}(\cH^{k_{\m}}\res E,p)\subseteq\bW_p$. With positive lower density, \cref{thm:c1_rec} shows that $E$ is $(k,\G)$-rectifiable of class $(1,0)$, and~\cref{thm:higher_order} concludes.
\end{proof}

\subsection{Sets of finite perimeter}\label{sec:finite_perimeter}

Throughout this subsection we assume $N\ge2$, so that hypersurfaces have
dimension $k=N-1\ge1$. In arbitrary Carnot groups, Don, Le Donne, Moisala, and Vittone~\cite{DonLeDonneMoisalaVittone2022} cover reduced boundaries of finite-perimeter sets by countably many sets with a cone property, a potentially weaker conclusion than rectifiability by regular surfaces.
We write $|A|:=\cL^N(A)$ for Haar volume, $d_{cc}$ for the
Carnot--Carath\'eodory distance associated with $X_1,\ldots,X_\ell$, and
$B^{cc}_r(x)$ for its open balls, and we follow~\cite{AKLD09}*{Definitions~4.11 and~4.15}.

\begin{definition}[Locally finite perimeter]
A measurable set $E\subset\G$ has
\emph{locally finite perimeter} if $X_1\chi_E,\ldots,X_\ell\chi_E$ are Radon
measures. Its perimeter measure $P_H(E;\cdot)$ is the total variation of
$D_H\chi_E:=(X_1\chi_E,\ldots,X_\ell\chi_E)$. Its \emph{reduced boundary}
$\mathcal FE$ is the set of $x$ such that $P_H(E;B^{cc}_r(x))>0$ for every
$r>0$, the limit
$\nu_E(x):=\lim_{r\downarrow0}D_H\chi_E(B^{cc}_r(x))/P_H(E;B^{cc}_r(x))$
exists, and $|\nu_E(x)|=1$.
\end{definition}

Changing $E$ on a Haar-null set changes none of these objects, so we may
assume that $E$ is Borel. The space $(\G,d_{cc},\cL^N)$ is complete, doubling,
and supports a weak $(1,1)$-Poincar\'e inequality~\cite{Jerison1986}, and
$P_H$ is its perimeter in the sense of~\cite{Miranda2003}. Indeed, for locally
Lipschitz $u$, $|\nabla_Hu|\le\operatorname{lip}u$ a.e.\ by Pansu's theorem,
with equality for smooth $u$, and smooth functions approximate in
perimeter~\cite{FSSC96}. Hence the results of~\cites{Amb01,Amb02} apply.

\begin{lemma}\label{lem:perimeter_facts}
Let $E\subset\G$ have locally finite perimeter. Then $\mathcal FE$ is Borel,
$\cH^{Q-1}\res\mathcal FE$ is locally finite, and for $\cH^{Q-1}$-a.e.\
$x\in\mathcal FE$ there is $t_x>0$ such that
\begin{equation}\label{eq:fp_densities}
 \begin{gathered}
 c\,t^{Q-1}\le\cH^{Q-1}(\mathcal FE\cap B_t(x))\le C\,t^{Q-1},\\
 \min\{|E\cap B_t(x)|,|B_t(x)\setminus E|\}\ge c\,t^Q,
 \end{gathered}
 \qquad 0<t<t_x,
\end{equation}
with $c,C>0$ depending only on $\G$ and $d$.
\end{lemma}

\begin{proof}
Since $r\mapsto P_H(E;B^{cc}_r(x))$ and $r\mapsto D_H\chi_E(B^{cc}_r(x))$ are
left-continuous and, for fixed $r$, Borel in $x$ (as differences of lower
semicontinuous functions), the limit defining $\nu_E$ may be taken along
rational radii, and $\mathcal FE$ is Borel.

By~\cite{AKLD09}*{Theorem~4.16 and~(4.17)}, which state the results
of~\cite{Amb01} in Carnot groups for sets of locally finite perimeter,
$P_H(E;\cdot)$ is concentrated on $\mathcal FE$ and has two-sided
$(Q-1)$-dimensional density bounds on $d_{cc}$-balls at $P_H(E;\cdot)$-a.e.\
point, with constants depending only on $\G$. By~\cite{AKLD09}*{(4.18)} it is
comparable to $\mathcal S^{Q-1}\res\mathcal FE$, built from $d_{cc}$ with the
normalization $\omega_{Q-1}$, hence to $\cH^{Q-1}\res\mathcal FE$, since
$\cH^{Q-1}\le\mathcal S^{Q-1}\le2^{Q-1}\cH^{Q-1}$ for the same distance and
normalization.

The volume bounds come from~\cite{Amb02}*{Theorem~5.4}, which concerns sets
of finite perimeter, so we truncate. For $x_0\in\G$ and $\rho'>0$ put
$B_\rho:=B^{cc}_\rho(x_0)$. The coarea formula~\cite{Miranda2003}*{Proposition~4.2}
for $d_{cc}(x_0,\cdot)$ gives $P_H(B_\rho;B_{\rho'})<\infty$ for almost every
$\rho<\rho'$, and then submodularity~\cite{Miranda2003}*{Proposition~4.7} gives
\[
 P_H(E\cap B_\rho;\G)=P_H(E\cap B_\rho;B_{\rho'})
 \le P_H(E;B_{\rho'})+P_H(B_\rho;B_{\rho'})<\infty.
\]
Inside $B_\rho$ the truncation has the same $D_H\chi$, reduced boundary, and
small-scale volume densities as $E$, so~\cite{Amb02}*{Theorem~5.4} gives
positive lower volume density of $E$ and of its complement at
$P_H(E;\cdot)$-a.e.\ point, with constants depending only on the doubling and
Poincar\'e constants. Finally, $d_{cc}$ and $d$ are bi-Lipschitz equivalent,
and all bounds transfer to the balls of $d$ and to our normalization
of~$\cH^{Q-1}$.
\end{proof}

In codimension one every vertical hyperplane subgroup is complemented by a
horizontal line (\cref{prop:vertical_splitting}), so~\cref{thm:main,thm:main_reverse}
apply with $k=N-1$ and $k_{\m}=Q-1$.

\begin{corollary}\label{cor:reduced_boundary}
Let $E\subset\G$ have locally finite perimeter, $m\in\mathbb N$, and
$\alpha\in[0,1]$. Then $\mathcal FE$ is $(N-1,\G)$-rectifiable of class
$(m,\alpha)$ with respect to $\cH^{Q-1}$ if and only if it satisfies
$\mathcal FE\cap B\in\mathfrak S_{m,\alpha}(N-1)$ for every ball
$B\subset\G$.
\end{corollary}

\begin{proof}
For every ball $B$, $\mathcal FE\cap B$ is Borel with finite
$\cH^{Q-1}$-measure by~\cref{lem:perimeter_facts}. If it belongs to
$\mathfrak S_{m,\alpha}(N-1)$, it is rectifiable of class $(m,\alpha)$
by~\cref{thm:main}: the positive lower-density hypothesis follows from~\cref{lem:perimeter_facts}, since sufficiently small balls centered in $B$ remain inside $B$. Countably many balls cover $\G$. Conversely, if
$\mathcal FE$ is rectifiable of class $(m,\alpha)$, so is each
$\mathcal FE\cap B$, and~\cref{thm:main_reverse} applies.
\end{proof}

The tubular condition can be verified from volume information alone.

\begin{lemma}
\label{lem:bulk_blowups}
Let $E\subset\G$ have locally finite perimeter, $m\in\mathbb N$, and
$\alpha\in[0,1]$. Suppose that for $\cH^{Q-1}$-a.e.\ $p\in\mathcal FE$ there
are an open neighborhood $\Omega_p$ of $p$ and $g_p\in\mathbf C^{m,\alpha}_\G(\Omega_p)$,
with $g_p(p)=0$ and $\nabla_Hg_p(p)\ne0$, such that
\begin{equation}\label{eq:bulk_rate}
 \bigl|(E\triangle\{g_p<0\})\cap B_r(p)\bigr|
 =o\bigl(r^{Q-1+m+\alpha}\bigr)
 \qquad\text{as }r\downarrow0.
\end{equation}
Then $\mathcal FE$ is $(N-1,\G)$-rectifiable of class $(m,\alpha)$
with respect to $\cH^{Q-1}$.
\end{lemma}

For instance, $g_p$ may be a polynomial of homogeneous degree at most $m$. For
$\alpha>0$ the general form does not reduce to the polynomial one, since
replacing $g_p$ by its Taylor polynomial changes $\{g_p<0\}\cap B_r(p)$ by a set
of volume $O(r^{Q-1+m+\alpha})$, while for $\alpha=0$ the two forms are equivalent.
These bounds follow from Taylor's theorem and the
estimate used in Step~4 of the proof of~\cref{prop:dirichlet_levels}.
The little-$o$ rate is needed in the argument: an $O(r^{Q-1+m+\alpha})$ rate
bounds the exceptional set in~\eqref{eq:main_assumption} only by a multiple of
$\lambda^{-1}r^{Q-1}$.

\begin{proof}
Fix $n\in\mathbb N$. Apply \cref{lem:compact_pieces} to
$\mu:=\cH^{Q-1}\res(\mathcal FE\cap B_n(e))$, with $\nu_1:=\cH^{Q-1}\res\mathcal FE$
(lower and upper bounds, exponent $Q-1$) and $\nu_2:=\cL^N\res E$,
$\nu_3:=\cL^N\res(\G\setminus E)$ (lower bounds, exponent $Q$),
and~\cref{lem:perimeter_facts}. Then $\mathcal FE\cap B_n(e)$ is covered, up to an
$\cH^{Q-1}$-null set, by compact sets $K$ for which there are $c',C',t_K>0$
such that~\eqref{eq:fp_densities} holds with $c=c'$ and $C=C'$ at every
$x\in K$ and every $t\in(0,t_K)$. To conclude the proof it is enough to show
that such a set $K$ is rectifiable of class $(m,\alpha)$.

Set $\nu:=\cH^{Q-1}\res\mathcal FE$. The bounds in~\cref{lem:perimeter_facts}
make $\nu$ locally asymptotically doubling. For such measures, closed balls
form a Vitali relation~\cite{Federer69}*{2.8.17}, and the density
theorem~\cite{Federer69}*{2.9.11}, together with the two-sided bounds to pass
from closed to open balls, gives
\[
 \frac{\nu(K\cap B_r(p))}{\nu(B_r(p))}\longrightarrow1
 \qquad\text{for $\nu$-a.e.\ }p\in K.
\]
Together with the lower bound for $\nu$, this proves
$\Theta_*^{Q-1}(K,p)>0$ for $\cH^{Q-1}$-a.e.\ $p\in K$.

Fix $p\in K$ satisfying~\eqref{eq:bulk_rate}, shrink $\Omega_p$ so that
$\nabla_Hg_p\ne0$ on it, and let $D_p:=\{q\in\Omega_p \ : \ g_p(q)<0\}$ and
$S_p:=g_p^{-1}(0)\cap\Omega_p$. Then $S_p$ is a $\G$-regular hypersurface of class
$(m,\alpha)$ through $p$, since in codimension one the kernel of $d_Hg_p(q)$
is complemented by a horizontal line (\cref{prop:vertical_splitting}). Fix $\lambda>0$ and set $\eta:=\lambda r^{m+\alpha}$,
$a:=\eta/10$, and
\[
 F_r:=K\cap B_r(p)\cap\bigl\{q \ : \ d(q,S_p)>\eta \bigr\}.
\]
If $m+\alpha=1$ and $\lambda\ge1$, then $F_r=\emptyset$ since $p\in S_p$;
otherwise $a<r/10$ for all small $r$.
Choose disjoint balls $B_a(q_i)$, with $q_i\in F_r$, whose fivefold dilates
cover $F_r$. For small $r$ these balls lie in $B_{2r}(p)\cap\Omega_p$ and
avoid $S_p$, and $5a<t_K$. Each ball is connected by radial dilation
paths, so $g_p$ has constant sign on it, and the volume bound
in~\eqref{eq:fp_densities} at $q_i\in K$ gives
$|(E\triangle D_p)\cap B_a(q_i)|\ge c'a^Q$. Since $F_r\subset\mathcal FE$,
the upper bound in~\eqref{eq:fp_densities} and disjointness yield
\[
 \cH^{Q-1}(F_r)
 \le C'\sum_i(5a)^{Q-1}
 \le \frac{C'5^{Q-1}}{c'}\,a^{-1}\bigl|(E\triangle D_p)\cap B_{2r}(p)\bigr|
 =o(r^{Q-1}).
\]
Thus $K$ satisfies the tubular condition of~\cref{def:special_class} at $p$
with $S=S_p$, for every $\lambda>0$, including the endpoint $\alpha=0$.
Hence $K\in\mathfrak S_{m,\alpha}(N-1)$. Since its positive lower density was proved above, applying~\cref{thm:main} with $k=N-1$
shows that $K$ is rectifiable of class $(m,\alpha)$ and concludes the proof.
\end{proof}

\begin{corollary}\label{cor:local_form}
Let $m\in\mathbb N$, $\alpha\in[0,1]$, and $\Omega\subset\G$ open, and let $E$
have locally finite perimeter in $\Omega$, with reduced boundary $\mathcal FE$
in $\Omega$. Then $\mathcal FE$ is $(N-1,\G)$-rectifiable of class $(m,\alpha)$
if and only if $\mathcal FE\cap B\in\mathfrak S_{m,\alpha}(N-1)$ for every ball
$B$ with $\overline B\subset\Omega$, and this holds if~\eqref{eq:bulk_rate}
holds at $\cH^{Q-1}$-a.e.\ $p\in\mathcal FE$.
\end{corollary}

\begin{proof}
The proofs of~\cref{lem:perimeter_facts,cor:reduced_boundary,lem:bulk_blowups}
apply verbatim with balls and compact pieces inside $\Omega$, taking
$\overline B^{cc}_{\rho'}(x_0)\subset\Omega$ in the truncation step.
\end{proof}

\begin{remark}\label{rmk:bulk_examples}
(i) For $m=1$ and $\alpha=0$, condition~\eqref{eq:bulk_rate} says that the
full blow-ups converge in $L^1_{\mathrm{loc}}$ to a vertical half-space. Whether
this holds at $\cH^{Q-1}$-a.e.\ point of $\mathcal FE$ is open in general
Carnot groups, but known in step two~\cite{FSSC03} and in some further
classes of groups~\cites{Marchi2014,LeDonneMoisala2021}. This identifies the hypothesis of~\cref{lem:bulk_blowups} at $(m,\alpha)=(1,0)$.

(ii) Let $D$ have locally finite perimeter in an open set $U$ and be
$K$-quasiminimal there, that is, $P_H(D;V)\le KP_H(F;V)$ whenever
$D\triangle F\Subset V\Subset U$, with $V$ open. Let $B$ be a $d_{cc}$-ball
with $\overline B\subset U$. Every $d_{cc}$-ball has finite perimeter by translation
and dilation of one such ball supplied by coarea. Submodularity therefore gives
finite perimeter for $D\cap B$. Given $V\Subset B$ open and a competitor $F$ for
$D\cap B$ with $(D\cap B)\triangle F\Subset V$, set $\widetilde F=F\cup(D\setminus B)$.
Then $D\triangle\widetilde F\Subset V$. Locality of perimeter gives
$P_H(D\cap B;V)=P_H(D;V)\le K P_H(\widetilde F;V)=K P_H(F;V)$,
as required by~\cite{KKLS2013}*{Definition~3.1}. For the
open representative of~\cite{KKLS2013}*{Theorem~4.2 and the discussion after
Remark~5.4}, $\partial D\cap U$ is the essential boundary, so
$\cH^{Q-1}((\partial D\setminus\mathcal FD)\cap U)=0$
by~\cite{Amb02}*{Theorem~5.3} and~\cite{AKLD09}*{Theorem~4.16}. Since
rectifiability ignores null sets, $\partial D\cap U$ is rectifiable of class
$(m,\alpha)$ if and only if $\mathcal FD\cap U$ is; by~\cref{cor:local_form},
this holds in particular if~\eqref{eq:bulk_rate} holds at
$\cH^{Q-1}$-a.e.\ point of $\mathcal FD\cap U$.
\end{remark}

Below, $W_H^{j,q}$ is the Sobolev space of weak horizontal derivatives of
order at most $j$.

\begin{proposition}
\label{prop:dirichlet_levels}
Let $\Omega\subset\G$ be open, let $g\in L^q_{\mathrm{loc}}(\Omega)$ for
some finite $q>\max\{Q,2\}$, and let
$u\in W^{1,2}_{H,\mathrm{loc}}(\Omega)$ be a weak solution of
$-\sum_{i=1}^{\ell}X_i^2u=g$ in $\Omega$. For instance, $u$ may locally
minimize the functional
\[
 \mathcal J_g(v;V):=\int_V\left(\frac12|\nabla_Hv|^2-gv\right)\, \mathrm dx,
 \qquad V\Subset\Omega,
\]
under compactly supported perturbations. For Lebesgue-a.e.\ $t\in\R$,
$E_t:=\{u<t\}$ has locally finite perimeter in $\Omega$, and its reduced
boundary in $\Omega$ is $(N-1,\G)$-rectifiable of class $(2,0)$.
If $\iota \le 2$, the covering hypersurfaces can be
chosen Euclidean $C^1$.
\end{proposition}

\begin{proof}
\proofstep{1 (Regularity)} Let $Q>2$. The fields $X_i$ are divergence-free,
so $X_i^*=-X_i$ and $u$ is a distributional solution. Since
$g\in L^q_{\mathrm{loc}}(\Omega)$, \cite{Folland1975}*{Theorem~6.1 and
Corollary~4.13} give $u\in W^{2,q}_{H,\mathrm{loc}}(\Omega)$,
and~\cite{Folland1975}*{Theorem~5.15} gives $u\in\Gamma_{2-Q/q}$ locally, with
$1<2-Q/q<2$. By the definition of $\Gamma_{1+\gamma}$ in~\cite{Folland1975},
$X_iu\in\Gamma_{1-Q/q}$ locally, so $u$ has a representative in
$\mathbf C_\G^1(\Omega)$. If $Q\le2$, the group is Euclidean of dimension at
most two, and classical interior estimates apply. Let
$H_u:=(X_iX_ju)_{i,j\le\ell}$. Each second-layer field $X_k$ is a combination of
commutators $[X_i,X_j]$, so the weak derivative $X_ku$ is the corresponding
combination of $X_iX_ju-X_jX_iu$, a fixed function in $L^q_{\mathrm{loc}}(\Omega)$.

\proofstep{2 (Pointwise second-order expansion)} Let $p$ be a common
$L^q$-Lebesgue point of these functions and $H_u$. By~\eqref{eq:Taylor_cond},
the values $u(p)$ and $X_iu(p)$ of the continuous representatives and the
Lebesgue values $X_ku(p)$ and $X_iX_ju(p)$, $i\le j\le\ell$, determine a
polynomial $T_p$ of homogeneous degree at most two. Then
$X_iX_jT_p\equiv(H_u(p))_{ij}$ also for $i>j$ since the relation
$[X_i,X_j]=\sum_kc^k_{ij}X_k$ holds for weak derivatives, so the Lebesgue
values satisfy $(H_u(p))_{ij}-(H_u(p))_{ji}=\sum_kc^k_{ij}X_ku(p)$, while
$X_kT_p\equiv X_ku(p)$ because $X_kT_p$ is constant when $d_k=2$.
The Poincar\'e inequality of~\cite{Jerison1986}, the bound
$|B_\rho|\approx\rho^Q$, and the usual telescoping argument give the scaled
Morrey inequality
\[
\operatorname{osc}_{B_\rho(x)}w\le C\rho\, \left(|B_{C\rho}(x)|^{-1}
\int_{B_{C\rho}(x)}|\nabla_Hw|^q \right)^{1/q} \qquad \text{for }q>Q.
\]
Apply this to the continuous functions $w=X_j(u-T_p)$, which vanish at $p$.
Since $p$ is an $L^q$-Lebesgue point of $H_u$, Morrey's inequality gives
$\sup_{B_{Cr}(p)}|\nabla_H(u-T_p)|=o(r)$. Integrating along horizontal curves
from $p$, of length at most $Cr$ and contained in $B_{Cr}(p)$, yields, after
enlarging the fixed constant $C$,
\[
 h_p(r):=\sup_{B_r(p)}|u-T_p|
 \le Cr^2\left(\frac{1}{|B_{Cr}(p)|}
 \int_{B_{Cr}(p)}|H_u(x)-H_u(p)|^q\, \mathrm dx\right)^{1/q}
 =o(r^2).
\]
\proofstep{3 (Coarea selection)} Let $Z_0$ be a Borel Haar-null set
containing all points that are not common Lebesgue points as above, and set $Z:=Z_0\cup\{\nabla_Hu=0\}$. The local coarea formula~\cite{Miranda2003}*{Proposition~4.2}
gives locally finite perimeter for almost every $E_t$ and, on a countable exhaustion of $\Omega$,
\[
\int_\R P_H(E_t;V\cap Z)\, \mathrm dt=\int_{V\cap Z}|\nabla_Hu|\,\mathrm dx=0.
\]
The perimeter--Hausdorff comparison therefore removes $Z$ from almost every reduced boundary. Fix a remaining point $p$ of such a boundary.

\proofstep{4 (Estimate)} By continuity $u(p)=t$, so $P_p:=T_p-t$ satisfies
$P_p(p)=0$ and $\nabla_HP_p(p)\ne0$. The sign mismatch lies in
$\{|P_p|\le h_p(r)\}\cap B_r(p)$. Choose horizontal $X$ with
$XP_p(p)\ne0$ and a vertical hyperplane $\mathbb W$ complementary to
$\exp(\R X)$. Choose a small product neighborhood in the coordinates
$x=pw\exp(sX)$ on which $XP_p$ keeps one sign and $|XP_p|$ is bounded below.
In these coordinates Haar measure is a constant multiple of the push-forward
of $|\cdot|_{\mathbb W}\otimes\cL^1$, since conjugation by $\exp(sX)$
preserves Haar measure on $\mathbb W$, as in the proof
of~\cref{lem:normal_form_norming}.
For small $r$, this neighborhood contains $B_r(p)$, and monotonicity in $s$
gives a slab of length at most $C_ph_p(r)$ on each fiber. Moreover,
the projection of $B_r(p)$ has Haar measure $O(r^{Q-1})$. Fubini then gives
\[
 |(E_t\triangle\{P_p<0\})\cap B_r(p)|
 \le C_ph_p(r)r^{Q-1}=o(r^{Q+1}).
\]
Apply~\cref{cor:local_form} with $m=2$, $\alpha=0$, and $g_p:=P_p$.
The Euclidean assertion, instead, follows from the brief argument in~\cref{sec:euclidean_regularity}.
\end{proof}

We include~\cref{prop:dirichlet_levels} as an instance in
which~\eqref{eq:bulk_rate} follows from an equation rather than as an assumption.
Its conclusion can also be obtained from the $\mathbf C^2_\G$ Lusin approximation
of~\cite{MR4637365}*{Theorem~1.1} and the coarea formula, since $u$ is twice
approximately differentiable a.e.; that route may allow weaker integrability
than $q>\max\{Q,2\}$, which we use for $\mathbf C^1_\G$ regularity (Step~1)
and the pointwise expansion (Step~2). Coarea arguments concern almost every
level of a fixed function, whereas~\cref{cor:reduced_boundary,lem:bulk_blowups}
apply to each individual set of locally finite perimeter.

\subsection{Euclidean regularity} \label{sec:euclidean_regularity}

Recall that $\iota$ is the step of $\G$. If $m\ge\iota$, every $k$-dimensional
$\G$-regular surface $S$ of class $(m,\alpha)$ is a $k$-dimensional $C^1$
submanifold of $\R^N$; hence, under the hypotheses of~\cref{thm:main}, $E$ is
covered, up to an $\cH^{k_{\m}}$-null set, by countably many such
submanifolds. Indeed, let $f\in[\mathbf C^{m,\alpha}_\G(\Omega)]^{N-k}$ define
$S$ locally. Since $d_l\le\iota\le m$ for every $l$, \cref{def:alpha_regular}
already requires continuous derivatives $X_lf$ along the whole adapted frame
$X_1,\ldots,X_N$. Mollifying as in the proof of~\cref{thm:taylor_estimate},
$f_\eps\to f$ and $X_lf_\eps\to X_lf$ locally uniformly. As each
$\partial_{x_j}$ is a combination of $X_1,\ldots,X_N$ with polynomial
coefficients, $f\in C^1(\Omega)$. Then $\nabla_Hf$ is the Euclidean Jacobian
$Df$ applied to $X_1,\ldots,X_\ell$, so surjectivity of $d_Hf$ gives
$\operatorname{rank}Df=N-k$, and the implicit function theorem applies. In
$\bH^1$, where $\iota=2$, the threshold is sharp. Indeed, there are
$\bH$-regular surfaces of class $1$ in $\bH^1$ with Euclidean Hausdorff
dimension $5/2$~\cite{KSC04}. We do not know whether $m\ge\iota$ is needed
in groups of higher step.

\appendix

\section{Measurable choice and compact pieces}\label{app:measurable}

We collect three standard measure-theoretic results. \Cref{lem:rational_radii}
is used in the proofs of~\cref{lem:trace_normalized_selection,thm:c1_rec},
\cref{lem:borel_choice} in the same two proofs, and~\cref{lem:compact_pieces}
in the proofs of~\cref{lem:trace_normalized_selection,thm:c1_rec,lem:bulk_blowups}.

\begin{lemma}\label{lem:rational_radii}
Let $\mu$ be a Borel measure on $\G$, let $p\in\G$ and $s,\gamma\ge0$, and let
$\varphi:\G\to[0,\infty]$ be Borel. For $L,r>0$ put
\[
 \beta_L(r):=r^{-s}\mu\bigl(B_r(p)\cap\{\varphi>Lr^\gamma\}\bigr).
\]
Then $\beta_{2^\gamma L}(r)\le2^s\beta_L(\varrho)$ whenever $r<\varrho<2r$. In
particular, for every $n\in\N$,
\[
 \sup_{0<r<1/(2n)}\beta_{2^\gamma L}(r)
 \le2^s\sup\bigl\{\beta_L(\varrho):\varrho\in\mathbb Q,\ 0<\varrho<1/n\bigr\}.
\]
\end{lemma}

\begin{proof}
If $r<\varrho<2r$, then $B_r(p)\subset B_\varrho(p)$,
$L\varrho^\gamma\le2^\gamma Lr^\gamma$, and $r^{-s}\le2^s\varrho^{-s}$. For the
second assertion choose a rational $\varrho\in(r,2r)$, so that $\varrho<1/n$.
\end{proof}

In the rest of this appendix, $\mu$ is a finite Borel measure on a Polish
space $X$, and $\mu$-measurable means measurable for the completion of $\mu$.

\begin{lemma}\label{lem:borel_choice}
Let $Y$ be a nonempty Polish space and let $\mathcal R\subset X\times Y$ be a Borel set
whose projection to $X$ has full $\mu$-measure. Then there is a Borel map
$y:X\to Y$ with $(x,y(x))\in\mathcal R$ for $\mu$-a.e.\ $x\in X$.
\end{lemma}

\begin{proof}
First, note that the projection $\mathcal P$ of $\mathcal R$ is analytic,
hence universally measurable~\cite{Kechris1995}*{Theorem~21.10}, so the
hypothesis is meaningful. The
Jankov--von Neumann uniformization theorem gives a universally measurable map
$y_0:\mathcal P\to Y$ with $(x,y_0(x))\in\mathcal R$;
see~\cite{Kechris1995}*{Theorems~18.1 and~29.9}. Extend $y_0$ by a constant
outside $\mathcal P$. By Lusin's theorem~\cite{Kechris1995}*{Theorem~17.12},
there are closed sets $F_h\subset X$ with $\mu(X\setminus F_h)<2^{-h}$ such
that $y_0|_{F_h}$ is continuous. Thus $y_0$ is Borel on the full-measure Borel
set $F:=\bigcup_hF_h$, and setting $y:=y_0$ on $F$ and $y$ constant off $F$
gives the required map.
\end{proof}

\begin{lemma}\label{lem:compact_pieces}
Let $X=\G$ and let $A\subset\G$ be $\mu$-measurable. Suppose we are given
\begin{enumerate}[label=\textup{(\alph*)}]
\item a $\mu$-measurable map $\Psi:A\to Y$ into a Polish space $Y$;
\item countably many sequences of $\mu$-measurable functions
      $g_{i,n}:A\to[0,\infty]$, $i\in I$, with $g_{i,n}\to0$ as
      $n\to\infty$ $\mu$-a.e.\ on $A$;
\item locally finite Borel measures $\nu_1,\ldots,\nu_{n_0}$ on $\G$, exponents
      $s_1,\ldots,s_{n_0}\ge0$, and sets $\mathcal J_-,\mathcal J_+\subset\{1,\ldots,n_0\}$,
      possibly overlapping, such that,
      for $\mu$-a.e.\ $x\in A$, $\Theta_*^{s_j}(\nu_j,x)>0$ for $j\in\mathcal J_-$ and
      $\Theta^{*s_j}(\nu_j,x)<\infty$ for $j\in\mathcal J_+$.
\end{enumerate}
If some of these data are not needed, take $Y$ to be a point, $I=\emptyset$, or
$n_0=0$. Then $A$ is covered, up to a $\mu$-null set,
by countably many compact sets $K\subset A$ such that $\Psi|_K$ is continuous,
$g_{i,n}\to0$ uniformly on $K$ for every $i\in I$, and, for some
$c_K,r_K>0$ and all $x\in K$ and $0<r<r_K$,
\[
    \nu_j(B_r(x))\ge c_Kr^{s_j}\quad \text{for } j\in \mathcal J_-,
    \qquad
    \nu_j(B_r(x))\le c_K^{-1}r^{s_j}\quad \text{for } j\in \mathcal J_+.
\]
\end{lemma}

\begin{proof}
After discarding a $\mu$-null set, we may assume that $A$ is Borel and that
the hypotheses hold at every point of $A$. For fixed $t>0$, the map
$x\mapsto\nu_j(B_t(x))$ is lower semicontinuous by Fatou's lemma. Hence, for
$h,k\in\mathbb N$, the set $D_{h,k}$ of points $x\in A$ with
$\nu_j(B_t(x))\ge t^{s_j}/h$ for $j\in\mathcal J_-$ and $\nu_j(B_t(x))\le ht^{s_j}$
for $j\in\mathcal J_+$, for every rational $t\in(0,1/k)$, is Borel, and these sets
cover $A$. Comparing an arbitrary $r\in(0,1/(2k))$ with rational radii in
$(r/2,r)$ and $(r,2r)$ shows that the same bounds hold on $D_{h,k}$ for all
such $r$, with $h$ replaced by $2^{s_*}h$, where
$s_*:=\max\{0,s_1,\ldots,s_{n_0}\}$. Thus the compact sets $K\subset D_{h,k}$ below satisfy the conclusion with $c_K:=(2^{s_*}h)^{-1}$ and
$r_K:=1/(2k)$.

Fix $h,k$ and $\eps>0$, and enumerate $I$ as a subset of $\mathbb N$.
Egorov's theorem~\cite{Folland1999}*{Theorem~2.33}
(applied to $\min\{g_{i,n},1\}$, which converges uniformly to $0$ on the same sets as $g_{i,n}$),
applied to the $i$-th sequence with error $2^{-i-2}\eps$, and Lusin's
theorem~\cite{Kechris1995}*{Theorem~17.12}, applied to $\Psi$ with error
$\eps/4$, give a $\mu$-measurable set $D'\subset D_{h,k}$ with
$\mu(D_{h,k}\setminus D')<\eps/2$ on which every $g_{i,n}\to0$ uniformly and
$\Psi$ is continuous. By inner regularity~\cite{Kechris1995}*{Theorem~17.11},
$D'$ contains a compact set $K$ with $\mu(D'\setminus K)<\eps/2$. Taking
$\eps=1/l$, $l\in\mathbb N$, and letting $h,k$ range over $\mathbb N$ gives
the countable family.
\end{proof}

\section*{Acknowledgments}
The authors thank Mario Santilli (University of L'Aquila) for stimulating discussions and useful feedback during the early stages of this project.

The mathematical results and arguments in this paper were developed by the authors over the course of the project, which began in 2024. Large language models were used at the end to assist with copyediting, presentation, and verification of the mathematical arguments; the authors take full responsibility for the final content.

%%%%%%%%%%%%%%%%%%%%%%%%%%%%%%%%%%%%%%%%%%%%%%%%%%%%%%%%%%%%%%%%%%%%%%%%%%%%%%%%%%%%%%%%%%%%%%%%
% Let bibliography pages end naturally rather than stretch reference spacing.
\raggedbottom
\bibliography{bibliography}

\end{document}